\documentclass[11pt]{article} 
\usepackage{latexsym} 
\usepackage{amsmath}
\usepackage{amssymb} 
\usepackage{amscd} 
\usepackage{graphicx}
\usepackage{ascmac}
\usepackage{mathrsfs}
\usepackage{mathrsfs}
\usepackage{pb-diagram}
\begin{document} 
%\Huge \bfseries \mathversion{bold}
%%%%%%%%%%%%%%%%%%%%%%%%%%%%%%%%% 
%%%%%%%%%%%%%%%%%%%%%%%%%%%%%%%%%
%%%%%%%%%%%%%%%%%%%%%%%%%%%%%%%%%
\newtheorem{Th}{Theorem}[section]
\newtheorem{Cor}{Corollary}[section]
\newtheorem{Prop}{Proposition}[section]
\newtheorem{Lem}{Lemma}[section]
\newtheorem{Def}{Definition}[section]
\newtheorem{Rem}{Remark}[section]
\newtheorem{Ex}{Example}[section]
\newtheorem{stw}{Proposition}[section]

%Definitions of bet ent

\newcommand{\bet}{\begin{Th}}
\newcommand{\ent}{\stepcounter{Cor}
   \stepcounter{Prop}\stepcounter{Lem}\stepcounter{Def}
   \stepcounter{Rem}\stepcounter{Ex}\end{Th}}

%Definitions of bec enc bep enp bel enl
%bef enf ber enr bee ene

\newcommand{\bec}{\begin{Cor}}
\newcommand{\enc}{\stepcounter{Th}
   \stepcounter{Prop}\stepcounter{Lem}\stepcounter{Def}
   \stepcounter{Rem}\stepcounter{Ex}\end{Cor}}
\newcommand{\bep}{\begin{Prop}}
\newcommand{\enp}{\stepcounter{Th}
   \stepcounter{Cor}\stepcounter{Lem}\stepcounter{Def}
   \stepcounter{Rem}\stepcounter{Ex}\end{Prop}}
\newcommand{\bel}{\begin{Lem}}
\newcommand{\enl}{\stepcounter{Th}
   \stepcounter{Cor}\stepcounter{Prop}\stepcounter{Def}
   \stepcounter{Rem}\stepcounter{Ex}\end{Lem}}
\newcommand{\bef}{\begin{Def}}
\newcommand{\enf}{\stepcounter{Th}
   \stepcounter{Cor}\stepcounter{Prop}\stepcounter{Lem}
   \stepcounter{Rem}\stepcounter{Ex}\end{Def}}
\newcommand{\ber}{\begin{Rem}}
\newcommand{\enr}{
   %\stepcounter{Rem} 
   \stepcounter{Th}\stepcounter{Cor}\stepcounter{Prop}
   \stepcounter{Lem}\stepcounter{Def}\stepcounter{Ex}\end{Rem}}
\newcommand{\bee}{\begin{Ex}}
\newcommand{\ene}{
 %\stepcounter{Ex}
   \stepcounter{Th}\stepcounter{Cor}\stepcounter{Prop}
   \stepcounter{Lem}\stepcounter{Def}\stepcounter{Rem}\end{Ex}}
\newcommand{\Proof}{\noindent{\it Proof\,}:\ }

%%%%%%%%%%%%%%%%%%%%%%%%%%%%%%%%%%%%%%%%%
%Beginning of Local Definition
%Local definitions
\newcommand{\EE}{\mathbb{E}}
\newcommand{\QQ}{\mathbb{Q}}
\newcommand{\R}{\mathbb{R}}
\newcommand{\C}{\mathbb{C}}
\newcommand{\ZZ}{\mathbb{Z}}
\newcommand{\KK}{\mathbb{K}}
\newcommand{\NN}{\mathbf{N}}
\newcommand{\PP}{\mathbb{P}}
\newcommand{\HH}{\mathbb{H}}
\newcommand{\OO}{\mathbb{O}}
\newcommand{\uuu}{\boldsymbol{u}}
\newcommand{\xxx}{\boldsymbol{x}}
\newcommand{\aaa}{\boldsymbol{a}}
\newcommand{\bbb}{\boldsymbol{b}}
\newcommand{\AAA}{\mathbf{A}}
\newcommand{\BBB}{\mathbf{B}}
\newcommand{\ccc}{\boldsymbol{c}}
\newcommand{\iii}{\boldsymbol{i}}
\newcommand{\jjj}{\boldsymbol{j}}
\newcommand{\kkk}{\boldsymbol{k}}
\newcommand{\rrr}{\boldsymbol{r}}
\newcommand{\FFF}{\boldsymbol{F}}
\newcommand{\yyy}{\boldsymbol{y}}
\newcommand{\ppp}{\boldsymbol{p}}
\newcommand{\qqq}{\boldsymbol{q}}
\newcommand{\nnn}{\boldsymbol{n}}
\newcommand{\vvv}{\boldsymbol{v}}
\newcommand{\eee}{\boldsymbol{e}}
\newcommand{\fff}{\boldsymbol{f}}
\newcommand{\www}{\boldsymbol{w}}
\newcommand{\0}{\boldsymbol{0}}
\newcommand{\lon}{\longrightarrow}
\newcommand{\ga}{\gamma}
\newcommand{\pa}{\partial}
\newcommand{\QED}{\hfill $\Box$}
\newcommand{\id}{{\mbox {\rm id}}}
\newcommand{\Ker}{{\mbox {\rm Ker}}}
\newcommand{\grad}{{\mbox {\rm grad}}}
\newcommand{\ind}{{\mbox {\rm ind}}}
\newcommand{\Real}{{\mbox {\rm Re}}}
\newcommand{\Aut}{{\mbox {\rm Aut}}}
\newcommand{\rot}{{\mbox {\rm rot}}}
\newcommand{\diver}{{\mbox {\rm div}}}
\newcommand{\Gr}{{\mbox {\rm Gr}}}
\newcommand{\LG}{{\mbox {\rm LG}}}
\newcommand{\PGL}{{\mbox {\rm PGL}}}
\newcommand{\GL}{{\mbox {\rm GL}}}
\newcommand{\SL}{{\mbox {\rm SL}}}
\newcommand{\Diff}{{\mbox {\rm Diff}}}
\newcommand{\Symp}{{\mbox {\rm Symp}}}
\newcommand{\Ct}{{\mbox {\rm Ct}}}
\newcommand{\Uns}{{\mbox {\rm Uns}}}
\newcommand{\rank}{{\mbox {\rm rank}}}
\newcommand{\sign}{{\mbox {\rm sign}}}
\newcommand{\Spin}{{\mbox {\rm Spin}}}
\newcommand{\Sp}{{\mbox {\rm Sp}}}
\newcommand{\SO}{{\mbox {\rm SO}}}
\newcommand{\Int}{{\mbox {\rm Int}}}
\newcommand{\Hom}{{\mbox {\rm Hom}}}
\newcommand{\codim}{{\mbox {\rm codim}}}
\newcommand{\type}{{\mbox {\rm type}}}
\newcommand{\ord}{{\mbox {\rm ord}}}
\newcommand{\Iso}{{\mbox {\rm Iso}}}
\newcommand{\corank}{{\mbox {\rm corank}}}
\def\mod{{\mbox {\rm mod}}}
\newcommand{\pt}{{\mbox {\rm pt}}}
\newcommand{\enP}{\hfill $\Box$ \par\vspace{5truemm}}
\newcommand{\qed}{\hfill $\Box$ \par}
\newcommand{\no}{\noindent}
\newcommand{\spe}{\vspace{0.4truecm}}
\newcommand{\Tan}{{\mbox {\rm Tan}}}
\newcommand{\SVC}{\mbox{\rm{SVC}}}
\newcommand{\End}{\mathcal{E}}

\newcommand{\dint}[2]{{\displaystyle\int}_{{\hspace{-0.45truecm}}{#1}}^{#2}}

\newcommand{\DEF}{$\stackrel{\mathrm{def}}{\Longleftrightarrow}\ $}
\newcommand{\kanren}{{\noindent\normalsize \ding{"2B}\ }}

\newcommand{\spcodim}{{\mbox {\rm sp-cod}}}
\newcommand{\diffcodim}{{\mbox {\rm diff-cod}}}

\title{Prolongations of $(3, 6)$-distributions by singular curves}

\author{Goo Ishikawa and Yoshinori Machida}

\date{ }

\maketitle

\begin{abstract}
A subbundle of rank $3$ in the tangent bundle over a $6$-dimensional manifold is called 
a $(3, 6)$-distribution if its local sections generate the whole tangent bundle by 
taking their Lie brackets once. An integral curve of a distribution, whose velocity vectors belong to the 
distribution, can be a singular curve or an abnormal extremal 
in the sense of geometric control theory. 
In this paper, given a $(3, 6)$-distribution, we prolong it, using the data of singular curves, to 
a $(3,5,7,8)$-distribution, to a $(3, 5, 7, 8, 9)$-distribution 
which possesses additional pseudo-product structure respectively. 
Regarding also another prolongation to a $(4, 6, 8)$-distribution, 
we show the equivalence of the classification problems of those four classes of distributions 
obtained from $(3, 6)$-distributions, 
generalising the correspondences of those in $B_3$-${\mathrm{SO}}(3,4)$-homogeneous models. 
\end{abstract}

\renewcommand{\thefootnote}{\fnsymbol{footnote}}

%\footnotetext{\scriptsize
%This work was supported by Kakenhi 24K06700. 
%}

\section{Introduction.}

Let $X$ be a $6$-dimensional manifold and 
$D \subset TX$ a distribution, i.e., a vector subbundle of the tangent bundle $TX$, of rank $3$ on $X$. 
Then $D$ is called a {\it $(3,6)$-distribution} if it has the growth vector $(3,6)$, namely, if 
$\rank(D^{(2)}) = 6$, i.e., $D^{(2)} = TX$, 
where $D^{(2)} := D + [D, D] = [D, D]$ is the derived system of $D$ by Lie bracket $[ \ , \ ]$ of local sections. 
The condition that $D$ is a $(3,6)$-distribution is equivalent to that 
there exists a local frame $\xi_1, \xi_2, \xi_3$ of $D$ such that 
$\xi_1, \xi_2, \xi_3, 
\xi_4 := [\xi_1, \xi_2], \xi_5 := [\xi_1, \xi_3], {\mbox{\rm and\  }} \xi_6 := [\xi_2, \xi_3]$
form a local frame of the whole tangent bundle $TX$. 
Note that $(3, 6)$-distributions form a generic class in the space of distributions of rank $3$ in a $6$-dimensional manifold, 
i.e., any small perturbation of a $(3, 6)$-distribution is a $(3, 6)$-distribution again. 
The typical example of $(3, 6)$-distribution is given by the 
$B_3$-${\mathrm{SO}}(3,4)$-homogeneous model (see \S \ref{B3-models}). 

In this paper we analyse $(3, 6)$-distributions using the notion of singular curves, 
or abnormal extremals, which arises from geometric control theory (see \cite{LS, Montgomery, AS}, 
see also \S\ref{Singular curves distributions},\S\ref{Singular curves $(3, 6)$-distributions}). 
By the prolongation of $(3, 6)$-distributions using their singular curves, 
we obtain $(3, 5, 7, 8)$-distributions with pseudo-product structure 
in the sense of Tanaka \cite{Tanaka, Yamaguchi}. 
Moreover the prolonged $(3, 5, 7, 8)$-distribution is re-prolonged, using singular curves again, 
to a $(3, 5, 7, 8, 9)$-distribution with pseudo-product structure. 
Furthermore we consider another prolongations to $(4, 6, 8)$-distributions from $(3, 6)$-distributions. 
The four classes of distributions thus obtained are prescribed by appropriate 
associated pseudo-product structures. 

Then we prove that the following four classes of distributions 
correspond bijectively and naturally to each other: 

\bet
\label{main-theorem}
There are bijective correspondences, compatible with prolongation procedures, 
between local isomorphism classes of $(3, 6)$-distributions, 
those of $(3, 5, 7, 8)$-distributions with pseudo-product structure of type $B_3(2,3)$, 
those of $(3, 5, 7, 8, 9)$-distributions with pseudo-product structure of type $B_3(1,2,3)$. 
Moreover there are also bijective correspondences, compatible with prolongation procedures, 
between local isomorphism classes of $(3, 6)$-distributions, 
those of $(4, 6, 8)$-distributions with pseudo-product structure of type $B_3(1,3)$ in the generalised sense 
and those of $(3, 5, 7, 8, 9)$-distributions with pseudo-product structure of type $B_3(1,2,3)$. 
\ent

The constructions performed in our paper are compared with the case of $B_3$, i.e., 
$\SO(3, 4)$-homogeneous spaces consisting 
of flags by null spaces in the space $\R^{3,4}$ with indefinite conformal metric of type $(3, 4)$ 
described in \S \ref{B3-models} and they are represented symbolically using the Dynkin diagrams: 
\[
\begin{diagram}
\node{}
\node{\bullet\ \vrule width 15pt height 3pt depth -2.3pt\ \bullet\Longrightarrow\bullet} 
\arrow{sw}\arrow{se} 
\node{}
\\
\node{\bullet\ \vrule width 15pt height 3pt depth -2.3pt\ \circ\Longrightarrow\bullet}\arrow{ne}\arrow{se} %\arrow{se} 
\node{}
\node{\circ\ \vrule width 15pt height 3pt depth -2.3pt\ \bullet\Longrightarrow\bullet}\arrow{nw}\arrow{sw} %\arrow{se} 
\\
\node{}
\node{\circ\ \vrule width 15pt height 3pt depth -2.3pt\ \circ\Longrightarrow\bullet}\arrow{ne}\arrow{nw}
\node{}
  \end{diagram}
\]
\begin{center}
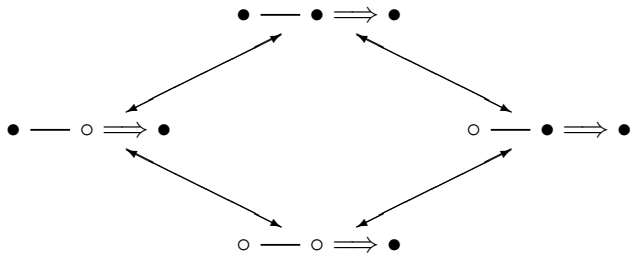
{\small Fig.1. Dynkin diagrams for prolongations-reductions from $(3,6)$-distributions. }
\end{center}

The bottom diagram $B_3(3)$ of the above diagram 
is associated with $(3, 6)$-distributions, the middle right $B_3(2, 3)$ 
with $(3, 5, 7, 8)$-distributions, the top $B_3(1, 2, 3)$ with $(3, 5, 7, 8, 9)$-distributions and the middle left  
$B_3(1, 3)$ with $(4, 6, 8)$-distributions respectively, 
together with appropriate pseudo-product structures in a generalised sense. 
Theorem \ref{main-theorem} is proved regarding these correspondences starting from general $(3, 6)$-distributions. 

In the process of prolongations via left routes from a $(3, 6)$-distribution in Figure 1, 
we observe that there arises the {\it generalised Legendre transformation} 
of $(3, 5, 7, 8, 9)$-distributions on the top. 
We remark the associated pseudo-product structures on the prolonged $(3, 5, 7, 8, 9)$-distributions 
by the two, right and left, routs of prolongations may differ in general but coincide in examples containing the flat case at least (see Remark \ref{generalised-Legendre-transformation} and Example \ref{example-4}). 

In \S \ref{Singular curves distributions} we recall the notions of singular characteristic (abnormal bi-extremals) 
and of singular curves (abnormal extremals) of distributions. 
We study on singular curves of $(3, 6)$-distributions in \S \ref{Singular curves $(3, 6)$-distributions}. 
We obtain the prolongations, which have the small growth vector $(3, 5, 7, 8)$ always, 
of $(3, 6)$-distributions by using their singular curves in \S \ref{Prolongation of $(3, 6)$-distributions by singular curves}, 
and investigate singular curves of prolonged $(3, 5, 7, 8)$-distributions of $(3, 6)$-distributions. 
Then we introduce naturally the class of $(3, 5, 7, 8)$-distributions with pseudo-product structure of type $B_3(2,3)$. 
In \S \ref{Iterated-prolongation of a $(3, 6)$-distribution} we observe that 
the iterated-prolongation by singular curves of the prolonged $(3, 5, 7, 8)$-distribution from 
a $(3, 6)$-distribution has the small growth vector vector $(3, 5, 7, 8, 9)$ always (Proposition \ref{iterated-prolongation}), 
and we introduce the class of $(3, 5, 7, 8, 9)$-distributions 
with pseudo-product structure of type $B_3(1,2,3)$. 
In \S \ref{Another-prolongation-of-(3, 6)-distribution}, we study another prolongations 
of $(3, 6)$-distributions to $(4, 6, 8)$-distributions, 
and introduce the class of $(4, 6, 8)$-distributions with pseudo-product structure of type $B_3(1,3)$ in a generalised sense. 
Note that the notion of pseudo-product structures on $(4, 6, 8)$-distributions 
is used in a weak sense than that by Tanaka \cite{Tanaka, Yamaguchi} (see 
\S \ref{Another-prolongation-of-(3, 6)-distribution}). 
We provide another way of prolongations of a $(4, 6, 8)$-distribution with pseudo-product structure of type $B_3(1,3)$ 
to a $(3, 5, 7, 8, 9)$-distribution directly using singular curves and discuss the difference of 
the corresponding $(3, 5, 7, 8, 9)$-distributions to the $(4, 6, 8)$-distribution in \S \ref{another-iterated-prolongation}. 
Then we prove Theorem \ref{main-theorem} in \S \ref{Proof of the main theorem}. 
In \S \ref{B3-models}, as an appendix, we describe the $B_3$-models of $(3, 5, 7, 8, 9)$-distributions, of $(3, 5, 7, 8)$-distributions, of $(4, 6, 8)$-distributions and of $(3, 6)$-distributions 
to compare those special cases of $B_3$-models with the general cases studied 
from \S \ref{Singular curves $(3, 6)$-distributions} to \S \ref{another-iterated-prolongation} 
of this paper.  

As prior researches on $(3, 6)$-distributions, 
Bryant \cite{Bryant, Bryant2} shows that the maximal symmetry of $(3, 6)$-distributions is realised by 
$B_3 = \SO(3, 4)$ symmetry on the model $D_0 \subset TX_0$ defined on the space $X_0 = {\mathcal F}_2$ 
which consists of null-planes of $\R^{3,4}$ (see also \S \ref{B3-models} of this paper). 
He shows that, given a $(3, 6)$-distribution $D \subset TX$, there exists an indefinite conformal 
structure of type $(3, 3)$ on $X$ associated to $D$ such that the invariants of $D$ are 
given by Weyl's conformal curvature tensor $W$ and the vanishing of $W$ is equivalent to that 
$(X, D) \cong (X_0, D_0)$ (see \cite{Bryant2}). 
This is compared with the existence of a conformal structure of type 
$(2, 3)$ on $(2, 3, 5)$-distributions discovered by Nurowski \cite{Nurowski, LNS}. 
It would be interesting to relate the data on singular curves of $(3, 6)$-distributions to those on 
Bryant's conformal $(3, 3)$-structures. Note that, generalising Bryant's result, 
Doubrov and Zelenko \cite{DZ2} solve the equivalence problem of 
$(3, 6, \dots)$-distributions under a weak condition. 
See also \cite{Krynski} for the related study on $(3, 5, 6)$-distributions. 

The classification problem of distributions has very long history and yields fruitful results from Cartan theory, 
representation theory of semi-simple Lie algebras, Tanaka theory and parabolic geometry 
(see \cite{Cartan, Bryant, LNS, CS}). 
The classical group $B_3$ is located on an interesting position between exceptional Lie groups $G_2$ and 
$F_4$. Regarding Dynkin diagrams symbolically, we observe that 
the $G_2$-diagram is obtained by \lq\lq folding\rq\rq\, 
from the $B_3$-diagram, and is embedded in the $F_4$-diagram. 
Also we observe that generic projections of a $(3, 6)$-distribution to a $5$-dimensional space 
turns to be $(3, 5)$-distributions coming from $(2, 3, 5)$-distributions which are well-studied in the relation with 
$G_2$-geometry (see also \cite{IMT, IKY, IKTY}). For example, $(S^3\times S^3)/(\ZZ/2)$ with 
a $(3, 6)$-distribution of $B_3(3)$-model projects to $(S^3\times S^3)/S^1 \cong S^3\times S^2$ 
with a $(2,3,5)$-distribution of $G_2(1)$-model. 
Moreover we observe that $(3, 6)$-distributions are embedded in a $(8, 15)$-distribution 
of type $F_4$ (\cite{IM2}). For example, $(3, 6)$-distributions arise from 
$\OO\PP^2_0$, a hyperplane section of the Cayley projective plane which possesses a $(8, 15)$-distribution 
of $F_4(4)$-model, by taking the projectivisation of tangent spaces to $\OO\PP^2_0$. These observations 
around $(3, 6)$-distributions could be studied in detail in forthcoming papers. 

In this paper, we do not treat the classification problems of $(3, 6)$-distributions directly, 
however, we provide a new aspect of singular curves from Hamiltonian geometry, 
and we show the equivalence of 
the classification problem of $(3, 6)$-distributions and those of related classes of distributions endowed with 
additional geometric structures, by giving bijective correspondences between those classes. 

In general, given a distribution $D \subset TX$, we set 
${\mathcal D}^{(2)} := {\mathcal D} + [{\mathcal D}, {\mathcal D}]$, 
${\mathcal D}^{(3)} := {\mathcal D}^{(2)} + [{\mathcal D}, {\mathcal D}^{(2)}]$, 
${\mathcal D}^{(4)} := {\mathcal D}^{(3)} + [{\mathcal D}, {\mathcal D}^{(3)}]$, $\dots$, 
regarding the sheaf ${\mathcal D}$ of local sections to $D$. 
Then the {\it small growth vector} of $D$ at $x_0 \in X$ is defined by the vector 
$(\rank({\mathcal D}), \rank(\pa^{(2)}{\mathcal D}), \rank(\pa^{(3)}{\mathcal D}), \dots)$ evaluated at $x_0$. 

In this paper, all manifolds and maps are supposed to be of class $C^\infty$ unless otherwise stated. 

\section{Singular curves of distributions.}
\label{Singular curves distributions}

First we recall the notions of singular curves and singular characteristics, or, abnormal extremals 
and abnormal bi-extremals (\cite{LS, Montgomery, AS}). 

Let $D \subset TX$ be any distribution on a manifold $X$ and $x_0 \in X$. 
Let ${\mathscr C}_{x_0}$ the space of $D$-integral curves, 
which consists of all absolutely continuous curves 
$\gamma : [a, b] \to X$ with $\gamma(a) = x_0, \gamma'(t) \in D_{\gamma(t)}$ for almost every $t \in [a, b]$. 
The endpoint mapping 
$\End : {\mathscr C}_{x_0} \to X$ is defined by $\End(\gamma) = \gamma(b)$ 
for $\gamma \in {\mathscr C}_{x_0}$. 
Then a {\it singular curve} or an {\it abnormal extremal} of $D$ means 
a critical point of $\End$ where the differential of $\End$ is not surjective. 
A singular curve of $D$ is called also a {\it $D$-singular curve}. 

Singular curves are characterised locally by a constrained Hamiltonian equation. 
Let $D$ is of rank $r$ and 
$\xi_1, \dots, \xi_r$ a local frame of $D$. 
Suppose $\xi_1, \dots, \xi_r$ are defied on $X$ for simplicity. 
We define the Hamiltonian function on the cotangent bundle $T^*X$ of the vector field $\xi_i$ by 
$H_{\xi_i}(x, p) = \langle p, \xi_i(x)\rangle$ for $(x, p) \in T^*X$ and
the Hamiltonian function of $D$ by 
$$
H(x, p; u) := u_1 H_{\xi_1}(x, p) + \cdots + u_r H_{\xi_r}(x, p), 
$$
where $u_1, \dots, u_r$ are control parameters. 
Then consider the Hamiltonian equation on $(x, p, u) = (x(t), p(t), u(t))$, 
$$
\dot{x} = \frac{\partial H}{\partial p}(x, p), \quad \ \dot{p} = - \frac{\partial H}{\partial x}(x, p)  \quad \cdots\cdots\cdots (*)
$$
with constraints $H_{\xi_1}(x, p) = 0, \dots, H_{\xi_r}(x, p) = 0$ and $p(t) \not= 0$. 
The constraints mean that $p(t) \in D_{x(t)}^\perp \setminus \{ 0\}$. Here 
$D^\perp \subset T^*X$ denotes the subbundle which consists of co-vectors annihilating $D$. 

We say that $(x(t), p(t))$ is a {\it singular characteristic} or an {\it abnormal bi-extremal} 
if there exist $u_1(t), \dots, u_r(t)$ depending on $t$ 
such that $(x(t), p(t), u(t))$ satisfies the equation $(*)$ with the constraints. 
Then a curve $\gamma(t) = x(t)$ in $X$ 
is a $D$-singular curve if and only if there exists a singular characteristic 
$\Gamma(t) = (x(t), p(t))$ covering $x(t)$ with $p(t) \not= 0$ (see \cite{AS, LS}). 
Note that a $D$-singular curve satisfies 
\[
\dot{x} = u_1(t)\xi_1(x(t)) + \cdots + u_r(t)\xi_r(x(t)), 
\]
with the controls $u_1(t), \dots, u_r(t)$. 

We define the Hamiltonian $H_{\xi} : T^*X \to \R$ of a vector field $\xi$ over $X$ by 
$H_{\xi}(x, p) := \langle p, \xi(x)\rangle$ for $(x, p) \in T^*X$, and 
denote by $\overrightarrow{H_{\xi}}$ Hamiltonian vector field over $T^*X$ of 
the Hamiltonian $H_{\xi}$. 
Note that $\overrightarrow{H_{\xi}}$ covers $\xi$ 
over $X$ via the canonical projection from $T^*X$ to $X$. 

If $(x(t), p(t), u(t))$ satisfies the equation $(*)$, 
then it holds along the curve 
\[
\frac{d}{dt}H_{\xi}(x(t), p(t)) = u_1(t) H_{[\xi_1, \xi]}(x(t), p(t)) + \cdots + u_r(t) H_{[\xi_r, \xi]}(x(t), p(t)) 
 \quad \cdots\cdots\cdots (**), 
\]
(see \cite{AS, LS}). 
Now we define the key notion of this paper: 

\bef
\label{singular-velocity-cone}
{\rm 
The {\it singular velocity cone} ${\mbox{\rm{SVC}}}(D) \subset TX$ of $D$ is defined as the set of 
$(x, u) \in TX$ such that there exists a $C^\infty$ $D$-singular curve $\gamma : (\R, 0) \to X$ such that 
$\gamma(0) = x, \gamma'(0) = u$. 
Moreover, given a subbundle $B \subset D^\perp (\subset T^*X)$, ${\mbox{\rm{SVC}}}(D, B)$ denotes the set of 
$(x, u) \in TX$ such that there exists a $D$-singular characteristic $\Gamma : (\R, 0) \to B$ with $\pi\circ\Gamma(0) = x$ 
and $(\pi\circ\Gamma)'(0) = u$, where $\pi : T^*X \to X$ is the projection. 
}
\enf

We see that, if $(x, u) \in TX$ is realised by the initial data of a singular curve $\gamma(t)$ 
at $t = 0$ and if $c \in \R^\times$, then $(x, cu)$ is realised by the singular curve $\gamma(ct)$, 
and therefore that ${\mbox{\rm{SVC}}}(D) \subseteq D$ and ${\mbox{\rm{SVC}}}(D)$ is a cone, 
i.e., it is invariant under the $\R^\times$-action on $TX$. Similarly we see ${\mbox{\rm{SVC}}}(D, B)$ is a sub-cone 
of $\SVC(D)$. 

\section{Singular curves of $(3, 6)$-distributions}
\label{Singular curves $(3, 6)$-distributions}

Let $D \subset TX$ be a $(3, 6)$-distribution. We study $D$-singular curves explained in the previous section. 
Consider the Whitney sum $T^*X \oplus D$ 
of the cotangent bundle $T^*X$ and the bundle $D$ over $X$, and 
define $H : T^*X \oplus D \to \R$ by 
\[
H(x, p, u) := 
u_1 H_{\xi_1}(x, p) + u_2 H_{\xi_2}(x, p) + u_3 H_{\xi_3}(x, p), 
\]
where $\{ \xi_1, \xi_2, \xi_3\}$ is a frame of $D$. 
Recall that $H_{\xi_i}$ is the Hamiltonian function defined on $T^*X$ by 
$H_{\xi_i}(x, p) := \langle p, \xi_i(x)\rangle$ for $i = 1, 2, 3$, as 
is introduced in \S \ref{Singular curves distributions}. 

Then $D$-singular curves are characterised by the Hamiltonian equation
\[
\left\{ 
\begin{array}{rcl}
\dot{x} & = & 
{\displaystyle 
\frac{\pa H}{\pa p} \ = u_1\xi_1(x(t)) + u_2\xi_2(x(t)) + u_3\xi_3(x(t)), 
}
\vspace{0.2truecm}
\\
\dot{p} & = & 
{\displaystyle \frac{\pa H}{\pa p} \ = 
- u_1\frac{\pa H_{\xi_1}}{\pa x} - u_2\frac{\pa H_{\xi_2}}{\pa x} 
- u_3\frac{\pa H_{\xi_3}}{\pa x}, 
}
\end{array}
\right. 
\]
with the constraints 
\[
H_{\xi_1}(x(t), p(t)) = 0, H_{\xi_2}(x(t), p(t)) = 0, H_{\xi_3}(x(t), p(t)) = 0, \quad p(t) \not= 0. 
\]
Differentiating the both sides of the constraint $H_{\xi_1}(x(t), p(t)) = 0$ by $t$, along singular curves, 
we have by $(**)$ in \S \ref{Singular curves distributions}, 
\[
0 = \frac{d}{dt}H_{\xi_1} = u_1H_{[\xi_1, \xi_1]} + u_2H_{[\xi_1, \xi_2]} + u_3H_{[\xi_1, \xi_3]}
= u_2H_{\xi_4} + u_3H_{\xi_5}. 
\]
Similarly we have 
\[
0 = \frac{d}{dt}H_{\xi_2} = u_1H_{[\xi_2, \xi_1]} + u_2H_{[\xi_2, \xi_2]} + u_3H_{[\xi_2, \xi_3]}
= - u_1H_{\xi_4} + u_3H_{\xi_6}, 
\]
and 
\[
0 = \frac{d}{dt}H_{\xi_3} = u_1H_{[\xi_3, \xi_1]} + u_2H_{[\xi_3, \xi_2]} + u_3H_{[\xi_3, \xi_3]}
= - u_1H_{\xi_5} - u_2H_{\xi_6}. 
\]

Let $\alpha_1, \dots, \alpha_6$ be the dual frame of $T^*X$ to the frame $\xi_1, \dots, \xi_6$ of $TX$. 
We set $p(t) = \sum_{i=1}^6 \varphi_i(t)\alpha_i(t)$. 
Then
\[
H_{\xi_i}(x(t), p(t)) = \langle p(t), \xi_i(x(t))\rangle = \varphi_i(t). 
\]
Thus we have $\varphi_1(t) = 0, \varphi_2(t) = 0, \varphi_3(t) = 0$ and 
\[
u_2\varphi_4 + u_3\varphi_5 = 0, \ u_1\varphi_4 - u_3\varphi_6 = 0, \ u_1\varphi_5 + u_2\varphi_6 = 0, \ \cdots\cdots\cdots\  (\dag)
\]
with $(\varphi_4, \varphi_5, \varphi_6) \not= (0, 0, 0)$. 
The equalities $(\dag)$ are equivalent to that 
\[
\left(
\begin{array}{ccc}
u_2 & u_3 & 0
\\
u_1 & 0 & -u_3
\\
0 & u_1 & u_2
\end{array}
\right)
\left(
\begin{array}{ccc}
\varphi_4
\\
\varphi_5
\\
\varphi_6
\end{array}
\right)
= 
\left(
\begin{array}{ccc}
0
\\
0
\\
0
\end{array}
\right), 
\]
and to that 
\[
\left(
\begin{array}{ccc}
0 & \varphi_4 & \varphi_5
\\
\varphi_4 & 0 & -\varphi_6
\\
\varphi_5 & \varphi_6 & 0
\end{array}
\right)
\left(
\begin{array}{ccc}
u_1
\\
u_2
\\
u_3
\end{array}
\right)
= 
\left(
\begin{array}{ccc}
0
\\
0
\\
0
\end{array}
\right). 
\]
We set $A = \left(
\begin{array}{ccc}
0 & \varphi_4 & \varphi_5
\\
\varphi_4 & 0 & -\varphi_6
\\
\varphi_5 & \varphi_6 & 0
\end{array}
\right)$. 
Then, for given any non-zero $(\varphi_4(t), \varphi_5(t), \varphi_6(t))$, 
the rank of $A$ is equals to $2$, and 
the solutions 
$(u_1(t), u_2(t), u_3(t))$ 
must be functional multiples of 
$(\varphi_6(t), - \varphi_5(t), \varphi_4(t))$. 

We set the vector field $\Xi$ over $T^*X$ by 
\[
\Xi := H_{\xi_6}\overrightarrow{H_{\xi_1}} - H_{\xi_5}\overrightarrow{H_{\xi_2}} + 
H_{\xi_4}\overrightarrow{H_{\xi_3}}, 
\]
where $\overrightarrow{H_{\xi_i}}$ denotes the Hamiltonian vector field over $T^*X$ with the Hamiltonian 
$H_{\xi_i}$. Note that the vector field $\overrightarrow{H_{\xi_i}}$ over $T^*X$ 
is a lifting of the vector field $\xi_i$ over $X$. 

Then we have
\[
\begin{array}{rcl}
\Xi(H_{\xi_1}) & = & H_{\xi_6}\overrightarrow{H_{\xi_1}}(H_{\xi_1}) 
- H_{\xi_5}\overrightarrow{H_{\xi_2}}(H_{\xi_1})
+ H_{\xi_4}\overrightarrow{H_{\xi_3}}(H_{\xi_1})
\\
& = & - H_{\xi_5}H_{[\xi_2, \xi_1]} + H_{\xi_4}H_{[\xi_3, \xi_1]} 
= H_{\xi_5}H_{\xi_4} - H_{\xi_4}H_{\xi_5} = 0. 
\end{array}
\]
Similarly we have $\Xi(H_{\xi_2}) = 0$ and $\Xi(H_{\xi_3}) = 0$. 
Since $D^\perp$ is a transverse intersection of $\{ H_{\xi_i} = 0\}, i = 1,2,3$, 
we see that the vector filed $\Xi$ is tangent to $D^\perp$ in $T^*X$ and 
the integral curves of $\Xi$ in $D^\perp \setminus (\mbox{\rm zero-section})$ projects to singular curves for $D$. 
Thus we have 
\bel
\label{D-singular-curves}
{\rm (\cite{Montgomery} \S 6.9)}\ 
For any $x \in X$ and for any direction $[u] \in P(D_x) (\subset P(T_xX))$, 
there exists a unique $D$-singular curve $\gamma : (\R, 0) \to X$ with 
$[\gamma'(0)] = [u]$ up to parametrisation. In particular we have $\SVC(D) = D$. 
\enl

\section{Prolongation of $(3, 6)$-distributions by singular curves.}
\label{Prolongation of $(3, 6)$-distributions by singular curves}

We construct the prolongation $(Z, E)$ of a $(3, 6)$-distribution $(X, D)$ using singular curves as follows: 
Set $Z = PD = (D \setminus \mbox{\rm{zero-section}})/\R^\times$, 
the fibrewise projectivisation of the bundle $D$ over the base $X$. 
Then $Z$ is of dimension $8$ and the induced projection 
$\pi_X : Z \to X$, $\pi(x, [u]) = x$ has the projective plane $\R P^2$ as fibres. 
We define the distribution $E \subset TZ$ over $Z$ by 
\[
E_{(x, [u])} := \{ w \in T_{(x, [u])}Z \mid \pi_*(w) \in [u]\}, 
\]
for $(x, [u]) \in Z$, where $[u]$ is regarded as a line in $T_xX$ through the origin. 
Note that $E$ is of rank $3$. 

Here we introduce the following notion after Noboru Tanaka \cite{Tanaka}: 

\bef
{\rm 
\label{B_3-pseudo-product}
Let $(Z, E)$ be a distribution with small growth vector $(3, 5, 7, 8)$ around a point on an $8$-dimensional manifold $Z$. 
\\
(1) 
A {\it pseudo-product structure of type $B_3(2,3)$} on $E$ is a direct sum decomposition 
$E = E_1 \oplus E_2$ into integrable subbundles $E_1$ and $E_2$ of rank $1$ and $2$ respectively such that 
\[
[E_1, E_2] = E^{(2)}, \ [E_1, E^{(2)}] = E^{(3)}, \ [E_2, E^{(2)}] = E^{(2)}, \ [E_1, E^{(3)}] = E^{(3)}, \ [E_2, E^{(3)}] = TZ. 
\]
(2) 
We say $E$ has a {\it gradation algebra of type $B_3(2,3)$} if there exists a local frame $\zeta_1, \zeta_2, \zeta_3$ of $E$ 
such that 
${[\zeta_1, \zeta_2]} \equiv: \zeta_4$\ mod.$E$, $[\zeta_1, \zeta_3] \equiv: \zeta_5$\ mod.$E$, 
$[\zeta_2, \zeta_3] \in \langle \zeta_2, \zeta_3\rangle$, 
$[\zeta_1, \zeta_4] \equiv: \zeta_6, 
[\zeta_1, \zeta_5] \equiv: \zeta_7, [\zeta_2, \zeta_4] \equiv 0, 
[\zeta_2, \zeta_5] \equiv 0,  [\zeta_3, \zeta_4] \equiv 0, [\zeta_3, \zeta_5] \equiv 0$, mod.$E^{(2)}$,  
$[\zeta_1, \zeta_6] \equiv 0, 
[\zeta_1, \zeta_7] \equiv 0, 
[\zeta_2, \zeta_6] \equiv 0, 
[\zeta_2, \zeta_7] \equiv \zeta_8, 
[\zeta_3, \zeta_6] \equiv \zeta_8, 
[\zeta_3, \zeta_7] \equiv 0$, mod.$E^{(3)}$. 
}
\enf

See \S \ref{B3-models}(ii) for an example of pseudo-product structures and gradation algebras of type 
$B_3(2,3)$. 

For our situation, we set $E_2 := \Ker(\pi_{X*} : TZ \to TX)$, which is an integrable subbundle of $E$ of rank $2$. 
Then we have 

\bel
\label{(3,5,7,8)-lemma} 
The prolongation $(Z, E)$ of a $(3, 6)$-distribution $(M, D)$ has 
the small growth vector $(3, 5, 7, 8)$. 
Moreover there exists unique complementary subbundle $E_1 \subset E$ of rank $1$ 
to $E_2$ such that the decomposition $E = E_1 \oplus E_2$ gives a pseudo-product structure of type 
$B_3(2,3)$ on $E$. 
\enl

\ber
\label{Cauchy-characteristic}
{\rm 
The subbundle 
$E_2$ (resp. $E_1$) is characterised as the Cauchy characteristic of $E^{(2)}$ (resp. $E^{(3)}$) as well (see \cite{BCGGG} for instance). 
Then $(X, D)$ is recovered as the quotient of $(Z, E^{(2)})$ by $E_2$. The local quotient of 
$(Z, E^{(3)})$ by $E_1$ will be studied in a forthcoming paper. 
}
\enr

\noindent
{\it Proof of Lemma \ref{(3,5,7,8)-lemma}}: 
Take any $(x_0, [u_0]) \in Z$ and fix it. Without loss of generality, we suppose the first component $u_1$ 
of the vector $u_0 \in D_{x_0}$ with respect to the frame $\xi_1, \xi_2, \xi_3$ is non-zero. 
On the open set $U = \{(x, [u]) \in Z \mid u_1 \not= 0\}$, we set 
$z_2 = u_2/u_1$ and $z_3 = u_3/u_1$, the inhomogeneous coordinates of $\pi_X$-fibres. 
Take a complementary subbundle $E_1$ of rank $1$ to $E_2$ in $E$ with a generator
\[
\zeta_1 = \xi_1 + z_2\xi_2 + z_3\xi_3 + a\, \zeta_2 + b\, \zeta_3, 
\]
where 
$\zeta_2 = \frac{\pa}{\pa z_2},  \ \zeta_3 = \frac{\pa}{\pa z_3}$ and $a, b$ are some functions on $U$. 
Then we have that $\zeta_1, \zeta_2, \zeta_3$ generate $E$ with 
$E_1 = \langle\zeta_1 \rangle$ and $E_2 = \langle\zeta_2, \zeta_3\rangle$. 
We will see in the process of the proof that $a, b$ are determined uniquely 
such that $E_1, E_2$ satisfy the required conditions. 

Now we have 
${[\zeta_1, \zeta_2]} = - \xi_2 - \zeta_2(a)\zeta_2 - \zeta_2(b)\zeta_3 
\equiv -\xi_2 =: \zeta_4$\ mod.$E$, $[\zeta_1, \zeta_3] \equiv - \xi_3 =: \zeta_5$\ mod.$E$, 
$[\zeta_2, \zeta_3] = 0$ 
and that $E^{(2)} = \langle \zeta_1, \zeta_2, \zeta_3, \zeta_4, \zeta_5\rangle$, which is of rank $5$. 

Moreover we have $[\zeta_1, \zeta_4] \equiv -\xi_4 + z_3\xi_6 =: \zeta_6, 
[\zeta_1, \zeta_5] \equiv - \xi_5 - z_2\xi_6 =: \zeta_7, [\zeta_2, \zeta_4] \equiv 0, 
[\zeta_2, \zeta_5] \equiv 0,  [\zeta_3, \zeta_4] \equiv 0, [\zeta_3, \zeta_5] \equiv 0$, mod.$E^{(2)}$ and therefore $E^{(3)} = \langle \zeta_1, \zeta_2, \zeta_3, \zeta_4, \zeta_5, \zeta_6, \zeta_7\rangle$, 
which is of rank $7$. 
Further we have $[E_1, E^{(2)}] = E^{(3)}, [E_2, E^{(2)}] = E^{(2)}$. 
Note that the quotient line bundle $TZ/E^{(3)}$ is generated by the class of $\xi_6$. 
Now fix the functions $a, b$ of $x, z_2, z_3$ on $U$ so that
\[
\begin{array}{c}
{[\zeta_1, \zeta_6]} \equiv [\xi_1 + z_2\xi_2 + z_3\xi_3, - \xi_4 + z_3\xi_6] + b \xi_6 \equiv 0, 
\\ 
{[\zeta_1, \zeta_7]} \equiv [\xi_1 + z_2\xi_2 + z_3\xi_3, - \xi_5 - z_2\xi_6] - a \xi_6 \equiv 0, 
\end{array}
\]
modulo $E^{(3)}$. 
Note that functions $a, b$ are uniquely determined, 
since the sheaf of sections to the line bundle $TZ/E^{(3)}$ is generated by the class of $\xi_6$ over $U$.
Lastly we set $\zeta_8 := -\xi_6$.  
Then we obtain $[\zeta_2, \zeta_6] \equiv 0, [\zeta_3, \zeta_6] \equiv \zeta_8, 
[\zeta_2, \zeta_7] \equiv \zeta_8, [\zeta_3, \zeta_7] \equiv 0$, mod.$E^{(3)}$. 
Then we have $[E_2, E^{(3)}] = E^{(4)} = TZ$. 
The local existence of $E_1$ implies the global existence of $E_1$. 

On other coordinate neighbourhoods, the arguments are performed similarly. 
\QED

\

In general we remark

\bep
\label{(3,5,7,8)-equivalence}
Let $E \subset TZ$ be a $(3,5,7,8)$-distribution. Then 
$E$ has a pseudo-product structure $E= E_1\oplus E_2$ of type $B_3(2,3)$ if and only if 
the gradation algebra of $E$ is of type $B_3(2,3)$. 
\enp

\Proof
Let $E$ has a gradation algebra of type $B_3(2,3)$. Then by setting $E_1 := \langle \zeta_1\rangle, E_2 := 
\langle \zeta_2, \zeta_3\rangle$, we see $E$ has a pseudo-product structure $E = E_1\oplus E_2$ of type 
$B_3(2,3)$. 

Conversely suppose $E$ has a pseudo-product structure $E = E_1\oplus E_2$ of type $B_3(2,3)$. 
Take any generator $\zeta_1$ of $E_1$ and any system of generators $\zeta_2, \zeta_3$ of $E_2$. 
Then we see that the conditions of Definition \ref{B_3-pseudo-product}(1) are fulfilled except for that 
$[\zeta_2, \zeta_7] - [\zeta_3, \zeta_6] \in E^{(3)}$, which follows in fact by 
$\zeta_7 \equiv [\zeta_1, [\zeta_1, \zeta_3]], \zeta_6 \equiv [\zeta_1, [\zeta_1, \zeta_2]], 
{\mbox {\rm mod.}} E^{(2)}$, by $[\zeta_1, [\zeta_1, [\zeta_2, \zeta_3]]] \equiv 0$, mod.$E^{(3)}$, 
and by 
the following purely algebraic reason (Lemma \ref{algebraic}). 
\QED

\

\bel
\label{algebraic}
Let $\zeta_1, \zeta_2, \zeta_3$ be sections of a distribution $E$. Then the equality 
\[
[\zeta_2, [\zeta_1, [\zeta_1, \zeta_3]]] - [\zeta_3, [\zeta_1, [\zeta_1, \zeta_2]]] = 
[\zeta_1, [\zeta_1, [\zeta_2, \zeta_3]]]
\]
holds. 
\enl

\Proof Using the Jacobi identity twice, we have 
\[
\begin{array}{l}
[\zeta_2, [\zeta_1, [\zeta_1, \zeta_3]]] - [\zeta_3, [\zeta_1, [\zeta_1, \zeta_2]]] 
\\
= 
\{ - [\zeta_1, [[\zeta_1, \zeta_3], \zeta_2]] - [[\zeta_1, \zeta_3], [\zeta_1, \zeta_2]]\} 
- \{ - [\zeta_1, [[\zeta_1, \zeta_2], \zeta_3]] - [[\zeta_1, \zeta_2], [\zeta_3, \zeta_1]]\} 
\\
= 
[\zeta_1, [\zeta_2, [\zeta_1, \zeta_3]]] - [\zeta_1, [\zeta_3, [\zeta_1, \zeta_2]]]
= 
- [\zeta_1, [\zeta_2, [\zeta_3, \zeta_1]]] + \{ [\zeta_1, [\zeta_1, [\zeta_2, \zeta_3]]] + [\zeta_1, [\zeta_2, [\zeta_3, \zeta_1]]]\}
\\
= [\zeta_1, [\zeta_1, [\zeta_2, \zeta_3]]]. 
\end{array}
\]
\QED

\

Now we see that there exists the natural bijection between the isomorphism classes of $(3,5,7,8)$-distributions with 
pseudo-product structures of type $B_3(2,3)$ and those of $(3, 6)$-distributions: 

\bep
\label{B_3(2,3)-and-(3,6)}
Let $E \subset TZ$ be a germ of $(3,5,7,8)$-distribution with a pseudo-product 
structure $E = E_1 \oplus E_2$ of type $B_3(2,3)$ on an $8$-dimensional manifold $Z$. 
Then the reduction $(X', D')$ of $(Z, E)$ by the integrable subbundle $E_2$ is a 
$(3, 6)$-distribution. Moreover the prolongation $(Z', E')$ of the reduced $(X', D')$ is 
isomorphic, as a distribution with pseudo-product structure of type $B_3(2,3)$,  
to the original distribution $(Z, E)$. 
Conversely, if $(Z, E)$, $E = E_1\oplus E_2$, is a prolongation of a $(3, 6)$-distribution $(X, D)$, 
then the reduction $(X', D')$ of $(Z, E)$ by $E_2$ is isomorphic to $(X, D)$. 
\enp

\Proof
Let $E_1 = \langle \zeta_1 \rangle, E_2 = \langle \zeta_2, \zeta_3\rangle$  as in Definition \ref{B_3-pseudo-product} 
and Lemma \ref{(3,5,7,8)-lemma} (cf. Proposition \ref{(3,5,7,8)-equivalence}). Then 
$E^{(2)} = \langle \zeta_1, \zeta_2, \zeta_3, \zeta_4, \zeta_5\rangle$. 
Then, by Remark \ref{Cauchy-characteristic}, 
the distribution $E^{(2)} \subset TZ$ is reduced to the distribution $D' := E^{(2)}/E_2 \subset T(M')$ of 
rank $3$ on the $6$-dimensional local leaf space $M' = Z/E_2$ of the integrable subbundle $E_2 \subset TZ$. 
The distribution $D'$ is generated by reduced vector fields 
$\overline{\zeta}_1, \overline{\zeta}_4, \overline{\zeta}_5$. 
We observe that $[\zeta_1, \zeta_4] = \zeta_6, [\zeta_1, \zeta_5] = \zeta_7$ and 
\[
\begin{array}{c}
[\zeta_4, \zeta_5] \equiv [[\zeta_1, \zeta_2], [\zeta_1, \zeta_3]] \equiv 
- [\zeta_1, [\zeta_3, [\zeta_1, \zeta_2]]] - [\zeta_3, [\zeta_1, [\zeta_1, \zeta_2]]]
\\
\equiv - [\zeta_1, [\zeta_3, \zeta_4]] - [\zeta_3, [\zeta_1, \zeta_4]]
\equiv - [\zeta_3, \zeta_6] \equiv - \zeta_8, 
\end{array}
\]
modulo $E^{(2)}$. Therefore we see that $[E^{(2)}, E^{(2)}] = TZ$. This implies that $[D', D'] = TM'$ and that 
$D' \subset TM'$ is a $(3, 6)$-distribution. 
Moreover the prolongation of $(M', D')$ is given, on $Z' = P(D')$, locally by 
$E' = \langle \zeta_1 + z_2\zeta_4 + z_3\zeta_5, \frac{\pa}{\pa z_2}, \frac{\pa}{\pa z_3}\rangle$ and 
we have that $(Z', E')$ is isomorphic to the original $(Z, E)$. 

The converse statement is clear. 
\QED

\bee
\label{example-1}
{\rm 
Here we give examples of $(3, 6)$-distributions $(X, D)$ and their prolongations $(Z, E), E = E_1\oplus E_2$ 
on $Z = P(D)$. 

On $X = \R^6 = \{ (x_1, x_2, x_3, x_4, x_5, x_6)\}$, consider distribution $D = \langle \xi_1, \xi_2, \xi_2\rangle$, 
where 
\[
\xi_1 = \frac{\pa}{\pa x_1} + x_3\frac{\pa}{\pa x_4} + (x_2 + m(x_6))\frac{\pa}{\pa x_6}, 
\  \ 
\xi_2 = \frac{\pa}{\pa x_2} + x_3\frac{\pa}{\pa x_5}, 
\ \ 
\xi_3 = \frac{\pa}{\pa x_3}, 
\]
and $m(x_6)$ is any $C^\infty$ function, that is not locally constant, on the variable $x_6$. 
Then we have that $E_1 = \langle \zeta_1\rangle, E_2 = \langle \frac{\pa}{\pa z_2}, \frac{\pa}{\pa z_3}\rangle$, 
where $\zeta_1 = \xi_1 + z_2\xi_2 + z_3\xi_3 + a\frac{\pa}{\pa z_2} + b\frac{\pa}{\pa z_3}$, 
$a = 0, b = - z_3m'(x_6)$. 
\qed
}
\ene

Now we study singular curves of the prolongation $(Z, E)$ of $(X, D)$. 
Recall the pseudo-product structure $E = E_1 \oplus E_2$ of $E$ obtained in Lemma \ref{(3,5,7,8)-lemma}. 

\bep
\label{SVC(E)}
The singular velocity cone of $E$ is given by $\SVC(E) = (E_1 \oplus 0) \cup (0 \oplus E_2)$. 
In fact, for any $(x_0, z_0) \in Z$ and for any $[w] \in P(E_{(x_0, z_0)})$, 
there exists a germ $\gamma : (\R, 0) \to Z$ of $C^\infty$ singular curve of $E$ 
satisfying $\gamma(0) = (x_0, z_0)$ and $[\gamma'(0)] = [w]$ if and only if 
$[w] \in (E_1 \oplus 0) \cup (0 \oplus E_2)$. Furthermore we have
\\
{\rm (1)} 
If $[w] \in P(E_1 \oplus 0)$, then $\gamma(t)$ is unique up to parametrisation and 
$\gamma(t)$ has bi-extremal lift in $E^{(3)\perp}$ with unique fibre components up to non-zero scalar multiplication. 
\\
{\rm (2)} 
If $[w] \in P(0 \oplus E_2)$, 
then $\gamma$ is any $C^\infty$ curve contained on $\pi_X$-fibre $(\pi_X)^{-1}(x_0) \cong \R P^2$ which satisfies
$\gamma(0) = (x_0, z_0)$ and $[\gamma'(0)] = [w]$.  
\enp

We have observed, in the case (2) of Proposition \ref{SVC(E)}, that any curve on any $\pi_X$-fibre is a singular curve of $E$. 
We see that any $\pi_X$-fibre is an {\it $E$-singular surface} if we introduce the following term: 

\bef
{\rm
Let $E \subset TZ$ be a distribution on a manifold $Z$. Then a submanifold $Z'$ of $Z$ is called 
an {\it $E$-singular submanifold} if any $E$-integral curve on $Z'$ is an $E$-singular curve. 
}
\enf

\

\noindent
{\it Proof of Proposition \ref{SVC(E)}:} 
The Hamiltonian of $(Z, E)$ in the sense of \S \ref{Singular curves distributions} 
is given, using the generators $\zeta_1, \zeta_2 = \frac{\pa}{\pa z_2}, \zeta_3 = \frac{\pa}{\pa z_3}$ 
of $E$ as in the Proof of Lemma \ref{(3,5,7,8)-lemma}, 
by
\[
H(x, z; p, r; \mu, \lambda) = \mu H_{\zeta_1} + \lambda_2H_{\zeta_2} + \lambda_3H_{\zeta_3} 
= \mu H_{\xi_1+z_2\xi_2+z_3\xi_3 + a\zeta_2 + b\zeta_3} + \lambda_2r_2 + \lambda_3r_3, 
\]
where $x = (x_1, \dots, x_6), z = (z_2, z_3); p = (p_1, \dots, p_6), r = (r_2, r_3)$ are local coordinates of $T^*Z$ 
and $\mu, \lambda_2, \lambda_3$ are control parameters. 
The Hamiltonian equation is given by
\[
\dot{x} = \frac{\pa H}{\pa p}, \ \dot{z} = \frac{\pa H}{\pa r}; \ \ 
\dot{p} = - \frac{\pa H}{\pa x}, \ \dot{z} = - \frac{\pa H}{\pa z}, 
\]
with constraints $H_{\zeta_1} = 0, H_{\zeta_2} = 0, H_{\zeta_3} = 0$ and $(p(t), r(t)) \not= 0$, 
namely, it is given by
\[
\begin{array}{c}
\dot{x} = \mu(\xi_1 + z_2\xi_2 + z_3\xi_3 + a\zeta_2 + b\zeta_3), \ \dot{z}_2 = \lambda_2 + a, \ \dot{z}_3 = \lambda_3 + b; 
\vspace{0.2truecm}
\\
\dot{p} = - \mu\, \frac{\pa}{\pa x} \left(H_{\xi_1+z_2\xi_2+z_3\xi_3 + a\zeta_2 + b\zeta_3}\right), \ \dot{r}_2 = - \mu H_{\xi_2}, \ \dot{r}_3 = - \mu H_{\xi_3}, 
\end{array}
\]
with constraints $r_2 = 0, r_3 = 0, H_{\xi_1} + z_2 H_{\xi_2} + z_3 H_{\xi_3} = 0$ and $(p(t), r(t)) \not= 0$. 

(1) First we treat the non-tangential case to $\pi_X$-fibre, i.e., $\mu \not= 0$. Then, by the constraints $r_2 = 0, r_3 = 0$, 
we have that $\mu H_{\xi_2} = 0, \mu H_{\xi_3} = 0$, so from $H_{\xi_1} + z_2 H_{\xi_2} + z_3 H_{\xi_3} = 0$, we have 
$H_{\xi_2} = 0, H_{\xi_3} = 0$ and then $H_{\xi_1} = 0$. 
Differentiate both sides of $H_{\xi_2} = 0$ by $t$, we have 
\[
0 = \mu H_{[\xi_2, \zeta_1]} + \lambda_2 H_{[\xi_2, \zeta_2]} + \lambda_3 H_{[\xi_2, \zeta_3]} 
= -\mu H_{\xi_4 - z_3\xi_6}, 
\]
since $[\xi_2, \zeta_2] = 0, [\xi_2, \zeta_3] = 0, r_2 = 0, r_3 = 0$. Similarly, from $H_{\xi_3} = 0$, we have 
\[
0 = \mu H_{[\xi_3, \zeta_1]} + \lambda_2 H_{[\xi_3, \zeta_2]} + \lambda_3 H_{[\xi_3, \zeta_3]} 
= - \mu H_{\xi_5 + z_2\xi_6}. 
\]
Since $\mu \not= 0$, we have $H_{\xi_4 - z_3\xi_6} = 0, H_{\xi_5 + z_2\xi_6} = 0$. 
Thus we see that $(p(t), 0) \in E^{(3)\perp}$. Note that $\rank(E^{(3)\perp}) = 1$. 

Let $I$ denote the ideal generated by 
$H_{\zeta_2}, H_{\zeta_3}, H_{\xi_1}, H_{\xi_2}, H_{\xi_3}, H_{\xi_4 - z_3\xi_6}, H_{\xi_5 + z_2\xi_6}$ 
in the ring of $C^\infty$ function germs on $T^*Z$ at a point of $Z$ and the corresponding co-vector $(p_0, 0) = (p(t_0), 0)$ 
for some $t_0$. Note that $H_{\zeta_1} \in I$. 
Then $E^{(3)\perp}$ is the zero-locus of the ideal $I$. 

Consider the Hamiltonian vector field $\Theta = \overrightarrow{H}_{\zeta_1}$ with Hamiltonian $H_{\zeta_1}$. 
Then we see all of the differentials $\Theta H_{\zeta_2}, \Theta H_{\zeta_3}, \Theta H_{\xi_1}, \Theta H_{\xi_2}, 
\Theta H_{\xi_3}, \Theta H_{\xi_4 - z_3\xi_6}, \Theta H_{\xi_5 + z_2\xi_6}$ belong to the ideal $I$. 
This means that the vector field $\Theta$ is tangent to the submanifold $E^{(3)\perp}$ of $T^*Z$. 
Thus we see integral curves of $\Theta$ with initial value in $E^{(3)\perp}$ is a singular bi-extremal
which projects to $\zeta_1$. This shows (1). 

(2) Next we treat the tangential case to $\pi_X$-fibre, i.e., $\mu = 0$. Let $\gamma : (\R, 0) \to Z$ 
be a germ of singular curve with an abnormal extremal $\Gamma : (\R, 0) \to T^*Z$ 
with $\Gamma(0) = (x_0, z_0; p_0, 0) \in T^*Z, p_0 \not= 0$, $\gamma'(0) 
\in \Ker(\pi_{X*})$. Assume $0$ is in the closure of $\Lambda = \{ t \mid \mu(t) \not= 0\}$. 
Then, by (1), $\gamma$ must be the 
integral curve of $\zeta_1$ at least in $\Lambda$. 
Beside, functions $a, b$ are continuous, so are locally bounded in $Z$
and therefore $\gamma$ is not tangential to $\pi_X$-fibres. This leads a contradiction. 
Now we suppose $\mu = 0$ on an open interval on $t$-line through $0$. 
In this case the Hamiltonian system is reduced to 
\[
\dot{x} = 0, \ \dot{z}_2 = \lambda_2 + a, \ \dot{z}_3 = \lambda_3 + b; \ \dot{p} = 0, \ \dot{r}_2 = 0, \ \dot{r}_3 = 0, 
\]
with constraint $r_2 = 0, r_3 = 0, H_{\xi_1 + z_2\xi_2+  z_3\xi_3} = 0$. 
Then $x, p$ must be constant,  $x(t) = x_0$ and $p(t) = p_0$, and $r(t)$ is equal to zero. 
Thus, in this case, the bi-extremal is of form $(x_0, z(t), p_0, 0)$ and any projected singular curve is contained in a $\pi_X$-fibre. 
Let $z(t) = (z_2(t), z_3(t))$ to be any pair of $C^\infty$ functions and $p_0 \in E^{(2)\perp}$.  Choose controls $\mu, \lambda_2, \lambda_3$ as 
$\mu = 0, \lambda_2 = z_2 - a, \lambda_3 = z_3 - b$. Then the constraint 
$H_{\xi_1} + z_2 H_{\xi_2} + z_3 H_{\xi_3} = 0$ is satisfied by $H_{\xi_1} = 0, H_{\xi_2} = 0, H_{\xi_3} = 0$ 
because $\xi_1, \xi_2, \xi_3 \in E^{(2)\perp}$. 
Hence the constrained Hamiltonian equation is satisfied and thus 
we see that any $C^\infty$ curve on a $\pi_X$-fibre is a singular curve. 
This shows (2). 

Combined with arguments in (1)(2), we have Lemma \ref{SVC(E)}. 
\QED

\ber
{\rm 
\label{SVC(E)-remark}
Note in Proposition \ref{SVC(E)} (2) that, if the bi-extremal $\Gamma$ satisfies that 
$p_0 \in E^\perp \setminus E^{(2)\perp}$, then the corresponding $E$-singular curve $\gamma(t) = (x(t), z_2(t), z_3(t))$ satisfies 
$x(t) = x_0$ and $\langle p_0, \xi_1(x_0)\rangle + z_2(t)\langle p_0, \xi_2(x_0)\rangle + z_3(t)\langle p_0, \xi_3(x_0)\rangle = 0$, 
which means that $z(t)$ is restricted to be an affine line in $z_2z_3$-plane of $\pi_X^{-1}(x_0)$. 
}
\enr

\section{Iterated-prolongation of a $(3, 6)$-distribution.}
\label{Iterated-prolongation of a $(3, 6)$-distribution}

In the previous sections, we have obtained the prolongation $(Z, E)$, that has the small growth vector 
$(3, 5, 7, 8)$, for a given $(3, 6)$-distribution $(X, D)$. 
Now we will prolong $(Z, E)$ again, by using singular curves, with 
directions in the sub-distribution $0\oplus E_2 \subset E = E_1\oplus E_2$ of rank $2$ 
to analyse the original $(3, 6)$-distributions more. 

Let $W = P(0\oplus E_2) \cong P(E_2)$, which is of dimension $9$ with the projection $\pi_Z : W \to Z$. 
Define the distribution $F \subset TW$ on $W$, for any $(x, [u], [\ell]) \in P(E_2), (x, [u]) \in Z, [\ell] \in P((0\oplus E_2)_{(x, [u])})$, by 
\[
F_{(x, [u], [\ell])} := \{ [v] \in T_{(x, [u], [\ell])}P(E_2) \mid \pi_{Z*}([v]) \in E_1\times [\ell] \}. 
\]
Then we see $F$ is a subbundle of $TW$ of rank $3$. In fact we have that $F$ is generated by 
\[
\theta_1 := \frac{\pa}{\pa w}, \ \ \theta_2 := \zeta_2 + w\zeta_3, \ \ \theta_3 := 
\xi_1 + z_2\xi_2 + z_3\xi_3 + a\zeta_2 + b\zeta_3 + c\frac{\pa}{\pa w}, 
\]
over a coordinate neighbourhood $\R \times U$ of $W$, 
where $\zeta_2 = \frac{\pa}{\pa z_2}, \zeta_3 = \frac{\pa}{\pa z_3}$, 
$w$ is a local coordinate of the $\pi_Z$-fibre, $c$ is a function on $\R\times U$, 
and $a, b$ are pull-backs by $\pi_Z$ of the functions defined on $U$ appeared in $\zeta_1$ 
and fixed uniquely in the proof of Lemma \ref{(3,5,7,8)-lemma}. 

We will see later soon that the small growth vector of $F$ is given by $(3,5,7,8,9)$. 
Furthermore we will show that the system of generators $\theta_1, \theta_2, \theta_3$ 
satisfy more strict properties: 

\bef
{\rm 
\label{B_3-(1,2,3)}
Let $F \subset TW$ be a germ of $(3,5,7,8,9)$-distribution. 
A {\it pseudo-product structure of type} $B_3(1,2,3)$ on $F$ 
is a direct-sum decomposition $F = F_1\oplus F_2\oplus F_3$ into subbundles $F_1, F_2, F_3$ of rank $1$, 
up to replacing $F_3$ by a line sub-bundle of $F_1\oplus F_3$ which is linearly independent of $F_1$, 
such that 
there exist generators $\theta_1, \theta_2, \theta_3$ of $F_1, F_2, F_3$ respectively satisfying the followings: 
\[
\begin{array}{l}
[\theta_1, \theta_2] =: \theta_4, \ [\theta_1, \theta_3] \in \langle \theta_1, \theta_3\rangle, 
\ [\theta_2, \theta_3] =: \theta_5, \ {\mbox{\rm so that }} F^{(2)} = \langle \theta_i \mid 1 \leq i \leq 5\rangle, \ 
\rank(F^{(2)}) = 5, 
\\
{[\theta_1, \theta_4]} \in \langle \theta_1, \theta_2, \theta_4\rangle, \ 
{[\theta_1, \theta_5]} \equiv: \theta_6\ {\mbox{\rm mod.}}F^{(2)}, \ 
[\theta_2, \theta_4] \in \langle \theta_1, \theta_2, \theta_4\rangle, \ 
{[\theta_2, \theta_5]} \in \langle \theta_2, \theta_3, \theta_5\rangle, 
\\
{[\theta_3, \theta_4]} \equiv -\theta_6\ {\mbox{\rm mod.}}F^{(2)}, \ 
[\theta_3, \theta_5] \equiv: \theta_7\  {\mbox{\rm mod.}}F^{(2)}, \ 
{\mbox{\rm so that }} F^{(3)} = \langle \theta_i \mid 1 \leq i \leq 7\rangle, \rank(F^{(3)}) = 7, 
\\
{[\theta_1, \theta_6]} \in \langle \theta_1, \theta_2,  \theta_3,  \theta_4, \theta_5,  \theta_6\rangle,\ 
{[\theta_1, \theta_7]} \equiv: \theta_8\ {\mbox{\rm mod.}}F^{(3)}, \ 
[\theta_2, \theta_6] \in \langle \theta_1, \theta_2,  \theta_3,  \theta_4, \theta_5,  \theta_6\rangle,\ 
\\
{[\theta_2, \theta_7]} \in \langle \theta_1, \theta_2,  \theta_3,  \theta_4, \theta_5,  \theta_7\rangle,\ 
{[\theta_3, \theta_6]} \equiv: \theta_8\ {\mbox{\rm mod.}}F^{(3)}, \ 
[\theta_3, \theta_7] \in F^{(3)},  \ [\theta_3, \theta_8] \in F^{(4)}, 
\\
{\mbox{\rm so that }} F^{(4)} = \langle \theta_i \mid 1 \leq i \leq 8\rangle, \rank(F^{(4)}) = 8, 
\\
{[\theta_1, \theta_8]} \in F^{(4)}, \ 
[\theta_2, \theta_8] \equiv: \theta_9\ {\mbox{\rm mod.}}F^{(4)}, \ 
{\mbox{\rm so that }} F^{(5)} = \langle \theta_i \mid 1 \leq i \leq 9\rangle = TW. 
\end{array}
\]
}
\enf

\ber
{\rm 
Note that the condition does not depend of the choice of generators $\theta_1, \theta_2, \theta_3$ of $F_1, F_2, F_3$ respectively. 
Also note that, if $F$ has a pseudo-product structure of type $B_3(1,2,3)$, then 
$F$ has isomorphic gradation algebra to that of $B_3(1,2,3)$ (see \S \ref{B3-models}(i)). 
}
\enr

\ber
\label{reduction-condition}
{\rm
The conditions of Definition \ref{B_3-(1,2,3)} imply that the distributions 
$\langle \theta_1, \theta_3\rangle, \langle \theta_1, \theta_2, \theta_4\rangle$, and 
$\langle \theta_1, \theta_2, \theta_3, \theta_4\rangle$ on $W$ 
are reduced to distributions on the local $\theta_1$-leaf space 
of rank $1, 2$ and $3$, respectively. Moreover we see that 
$\langle \theta_1, \theta_2, \theta_4\rangle, \langle \theta_2, \theta_3, \theta_5\rangle$ 
and $\langle \theta_1, \theta_2, \theta_3, \theta_4, \theta_5\rangle = F^{(2)}$ on $W$ 
are reduced to a distribution on the local $\theta_2$-leaf space of rank $2, 2$ and $4$, respectively. 
}
\enr

\bep
\label{iterated-prolongation}
For any germ of $(3, 6)$-distribution $(X, D)$, $x_0 \in X$ and at any point on $(\pi_X\circ\pi_Z)^{-1}(x_0)$ in $W$, 
the iterated prolongation $(W, F)$ of $(X, D)$ 
has $(3, 5, 7, 8, 9)$ as the small growth vector and a pseudo-product structure
$F = F_1\oplus F_2 \oplus F_3$ 
of type $B_3(1,2,3)$ in the sense of Definition \ref{B_3-(1,2,3)}. 
\enp

\Proof
Since  
$[\theta_1, \theta_2] = \zeta_3 =: \theta_4, \ [\theta_1, \theta_3] = c_w\frac{\pa}{\pa w} \equiv 0$ mod.$F$, 
$\ [\theta_2, \theta_3] \equiv \xi_2 + w\xi_3 + (\zeta_2 + w\zeta_3)(a)\zeta_2 + (\zeta_2 + w\zeta_3)(b)\zeta_3  - c\zeta_3 =: \theta_5$, mod.$F$, 
we see that $F^{(2)} = \langle \theta_1, \theta_2, \theta_3, \theta_4,  \theta_5\rangle$, 
which is a subbundle in $TW$ of rank $5$. 

Moreover we have $[\theta_1, \theta_4] = 0,  \ 
[\theta_1, \theta_5] \equiv \xi_3 + \zeta_3(a)\zeta_2 + \zeta_3(b)\zeta_3 - a_w\zeta_3 
\equiv \xi_3 =: \theta_6$, mod.$F^{(2)}$, \ 
$[\theta_2, \theta_4] = 0, \ 
[\theta_2, \theta_5] \equiv 0, \ 
[\theta_3, \theta_4] \equiv -\xi_3 = -\theta_6, \ 
[\theta_3, \theta_5] \equiv \xi_4 - z_3\xi_6 + w(\xi_5 + z_2\xi_6) - (\zeta_2 + w\zeta_3)(a)\xi_2 - (\zeta_2 + w\zeta_3)(b)\xi_3 +  c\xi_3 =: \theta_7$, mod.$F^{(2)}$, 
and we see that $F^{(3)} = 
 \langle \theta_i \mid 1 \leq i \leq 7\rangle$, which is 
 a subbundle in $TW$ of rank $7$. Note that $F^{(3)} \supset \langle \zeta_2, \zeta_3, \xi_1, \xi_2, \xi_3\rangle$. 
 
Furthermore, we have 
$[\theta_1, \theta_6] \equiv 0$, mod.$\langle\theta_1, \theta_2, \theta_3, \theta_4, \theta_5, \theta_6\rangle$, \ 
$[\theta_1, \theta_7] \equiv \xi_5 + z_2\xi_6 - \zeta_3(a)\xi_3 - \zeta_3(b)\xi_3 + c_w\xi_3 \equiv \xi_5 + z_2\xi_6 =:\theta_8, \ 
[\theta_2, \theta_6] \equiv 0, \ 
[\theta_2, \theta_7] \equiv 0, \ 
[\theta_3, \theta_6] \equiv \xi_5 + z_2\xi_6 = \theta_8,$ mod.$F^{(3)}$. Moreover, 
$[\theta_3, \theta_7] = 
\left(
[\xi_1 + z_2\xi_2 + z_3\xi_3, \xi_4 - z_3\xi_6] - b\xi_6
\right)
+ w
\left(
[\xi_1 + z_2\xi_2 + z_3\xi_3, \xi_5 + z_2\xi_6] + a\xi_6
\right) + c(\xi_5 + z_2\xi_6) 
\equiv 0$, ${\mbox {mod.}}F^{(3)}$, if we take $c$ properly. 
Note that $\pi_Z^*E^{(3)} \subset \langle \theta_i \mid 1 \leq i \leq 8\rangle = F^{(3)} + \langle\theta_8\rangle$ (see \S \ref{Prolongation of $(3, 6)$-distributions by singular curves}). 
Then we have 
$F^{(4)} = \langle \theta_i \mid 1 \leq i \leq 8\rangle$, which is of rank $8$. 

Lastly, by 
$[\theta_1, \theta_8] \equiv 0, \ 
[\theta_2, \theta_8] \equiv \xi_6 =: \theta_9, \ 
[\theta_3, \theta_8] \equiv - [\xi_1 + z_2\xi_2 + z_3\xi_3, \xi_5 + z_2\xi_6] - a\xi_6 \equiv 0$, {\mbox {mod.}}$F^{(4)}$, 
since $\pi_Z^*E^{(3)} \subset F^{(4)}$, and 
we have that $F^{(5)} = \langle \theta_i \mid 1 \leq i \leq 9\rangle = TW$. 

By the above calculations, we have that the small growth vector of $F$ is given by 
$(3, 5, 7, 8, 9)$ and, by setting $F_i = \langle \theta_i\rangle, i = 1,2,3$, we see that 
$F = F_1\oplus F_2\oplus F_3$ gives a pseudo-product structure of type $B_3(1,2,3)$ as in Definition 5.1. 
\QED

\

\ber
{\rm 
Note that necessarily we have $[\theta_3, \theta_4] \equiv - [\theta_1, \theta_5]$, mod.$F^{(2)}$, and 
$[\theta_3, \theta_6] \equiv [\theta_1, \theta_7]$, mod.$F^{(3)}$, 
by the following Lemma \ref{algebraic-calculation}. 
}
\enr

\bel
\label{algebraic-calculation}
Let $\theta_1, \theta_2, \theta_3$ be a section of a distribution $F$. If 
$[\theta_1, \theta_3] \equiv 0$, {\mbox{\rm mod.}}$F$, then we have 
\[
{[\theta_3, [\theta_1, \theta_2]]} \equiv - [\theta_1, [\theta_2, \theta_3]], {\mbox{\rm mod.}}F^{(2)}, \quad 
[\theta_3, [\theta_1, [\theta_2, \theta_3]]] \equiv [\theta_1, [\theta_3, [\theta_2, \theta_3]]], {\mbox{\rm mod.}}F^{(3)}. 
\]
\enl
\vspace{-0.2truecm}
\Proof 
By Jacobi identity, we have 
${[\theta_3, [\theta_1, \theta_2]]} = - [\theta_1, [\theta_2, \theta_3]] - [\theta_2, [\theta_3, \theta_1]] 
\equiv - [\theta_1, [\theta_2, \theta_3]]$, mod.$F^{(2)}$, 
and 
$
[\theta_3, [\theta_1, [\theta_2, \theta_3]]] = 
- [\theta_1,[[\theta_2, \theta_3], \theta_3]] - [[\theta_2, \theta_3], [\theta_3, \theta_1]]
\equiv - [\theta_1,[[\theta_2, \theta_3], \theta_3]] = [\theta_1, [\theta_3, [\theta_2, \theta_3]]]$, mod.$F^{(3)}$. 
\QED

\bee
\label{example-2}
{\rm
Continued with Example \ref{example-1}, the prolongation $(W, F)$ 
has the decomposition in Proposition \ref{iterated-prolongation} given by 
$F = F_1\oplus F_2\oplus F_3 = \langle \theta_1\rangle\oplus \langle \theta_2\rangle\oplus\langle \theta_3\rangle$ with 
$\theta_1 = \frac{\pa}{\pa w}, \theta_2 = \frac{\pa}{\pa z_2} + w\frac{\pa}{\pa z_3}, \theta_3 = 
\xi_1 + z_2\xi_2 + z_3\xi_3 + a\frac{\pa}{\pa z_2} + b\frac{\pa}{\pa z_3} + c\frac{\pa}{\pa w}$
with $a = 0, b = - z_3m'(x_6), c = - w m'(x_6)$. \qed
}
\ene

Now we study singularities of $F$ in detail. For a sub-bundle $G \subset F^\perp \subset T^*M$, 
we denote by $\SVC(F, G)$ the set of velocity vectors of singular curves with singular characteristics 
in $G$. See \S \ref{Singular curves distributions} for the definition of singular characteristics 
(see Definition \ref{singular-velocity-cone}). 
Note that $\SVC(F, F^\perp) = \SVC(F)$. Then we have

\bep
\label{(1,2,3)-SVC}
The singular velocity cone $\SVC(F)$ of $F = \{ \lambda_1\theta_1(w) + \lambda_2\theta_2(w) + \lambda_3\theta_3(w) 
\mid w \in W, \lambda_1, \lambda_2, \lambda_3 \in \R\}$ is given by 
$\{ \lambda_2 = 0\} \cup \{ \lambda_3 = 0\} = (F_1\oplus F_3)\cup(F_1\oplus F_2)$. 
Moreover we have 
\[
\begin{array}{c}
\SVC(F, F^\perp\setminus F^{(2)\perp}) = F_1\oplus F_3, \quad 
\SVC(F, F^{(2)\perp}\setminus F^{(3)\perp}) = F_1\oplus F_2, 
\\
\SVC(F, F^{(3)\perp}\setminus F^{(4)\perp}) = F_2, \quad {\mbox{\rm and}} \quad
\SVC(F, F^{(4)\perp}\setminus F^{(5)\perp}) = F_1\oplus F_3. 
\end{array}
\]
In particular, the $F_1$ and $F_2$ are intrinsically determined, 
and $F_3$ is determined up to replacing by another line sub-bundle of $F_1\oplus F_3$ 
which is independent of $F_3$. 
\enp

\Proof
For the Hamiltonian 
$
H = \lambda_1 H_{\theta_1} + \lambda_2 H_{\theta_2} + \lambda_3 H_{\theta_3}, 
$
we solve the Hamiltonian equation 
\[
\dot{x} = \frac{\pa H}{\pa p}, \ \ \dot{z} = \frac{\pa H}{\pa q}, \ \ \dot{w} = \frac{\pa H}{\pa r}, \quad
\dot{p} = - \frac{\pa H}{\pa x}, \ \ \dot{q} = - \frac{\pa H}{\pa z}, \ \ \dot{r} = - \frac{\pa H}{\pa w}, 
\]
with constraints $H_{\theta_1} = 0, H_{\theta_2} = 0, H_{\theta_3} = 0$ 
in $T^*W \setminus ({\mbox{\rm zero-section}})$ 
to find singular characteristics which project to singular curves in $W$. 

Along singular characteristics, we have 
\[
\frac{d}{dt}H_{\theta_1} = \lambda_1 H_{[\theta_1, \theta_1]} + \lambda_2 H_{[\theta_2, \theta_1]} + \lambda_3 H_{[\theta_3, \theta_1]} = - \lambda_2 H_{\theta_4}, 
\]
using the bracket relations obtained in the proof of Proposition \ref{iterated-prolongation}. 
Similarly we have 
\[
\begin{array}{c}
{\displaystyle 
\frac{d}{dt}H_{\theta_2} = \lambda_1 H_{\theta_4} - \lambda_3 H_{\theta_5}, 
\ \ 
\frac{d}{dt}H_{\theta_3} = \lambda_2 H_{\theta_5}, 
\ \ 
\frac{d}{dt}H_{\theta_4} = - \lambda_3 H_{\theta_6}, 
\ \ 
\frac{d}{dt}H_{\theta_5} = \lambda_1 H_{\theta_6} + \lambda_3 H_{\theta_7}, 
}
\vspace{0.1truecm}
\\
{\displaystyle 
\frac{d}{dt}H_{\theta_6} = \lambda_3 H_{\theta_8}, 
\ \
\frac{d}{dt}H_{\theta_7} = \lambda_1 H_{\theta_8}, 
\ \
\frac{d}{dt}H_{\theta_8} = \lambda_2 H_{\theta_9}. 
}
\end{array}
\]

From the constraints, on the closure of the domain $(H_{\theta_4}, H_{\theta_5}) \not= (0, 0)$, 
we have that $\lambda_2 = 0$ and that $\lambda_1 : \lambda_3 = H_{\theta_5} : H_{\theta_4}$ on $F^\perp \setminus F^{(2)\perp}$. 
Consider the non-singular 
vector field $\Theta = H_{\theta_5}\overrightarrow{H}_{\theta_1} + H_{\theta_4}\overrightarrow{H}_{\theta_3}$ 
over $F^\perp \setminus F^{(2)\perp}$. 
Then we see $\Theta H_{\theta_i}$, $i = 1, 2, 3,$ belong to the ideal generated by $H_{\theta_1}, H_{\theta_2}, H_{\theta_3}$ 
in the algebra of $C^\infty$ functions on $T^*W$. Thus $\Theta$ is tangent to $F^\perp \setminus F^{(2)\perp}$ 
and integral curve-germs to $\Theta$ are contained in $F^\perp \setminus F^{(2)\perp}$, which are singular curves of $F$.  
Note that the directions of singular curves of this case form the subbundle $\langle \theta_1, \theta_3\rangle$ generated by 
$\theta_1, \theta_3$. 

Now suppose $H_{\theta_4} = 0, H_{\theta_5} = 0$ along a singular characteristic on an interval. 
Then we have 
$\lambda_3 H_{\theta_6} = 0, \lambda_1 H_{\theta_6} + \lambda_3 H_{\theta_7} = 0$. 

Suppose $H_{\theta_6} \not= 0$. Then $\lambda_3 = 0$ and then $\lambda_1 = 0$. 
We set $\Theta' = \overrightarrow{H}_{\theta_2}$. 
Then we see that all of $\Theta' H_{\theta_i}, 1 \leq i \leq 5$ belong the ideal generated by 
$H_{\theta_j}, 1 \leq j \leq 5$ and that the integral curves of $\theta_2$ are singular curves with 
singular bi-characeristics in $F^{(2)} \setminus F^{(3)}$. Thus $\lambda_2$-axis is the singular direction in this case. 

Suppose  $H_{\theta_6} = 0, H_{\theta_7} \not= 0$ on an interval. Then $\lambda_3 = 0$ on that interval. 
Then the vector field 
$\Theta'' = A\overrightarrow{H}_{\theta_1} + B\overrightarrow{H}_{\theta_2}$, where $A, B$ are arbitrary functions 
on $T^*W$ with $(A, B) \not= (0, 0)$, satisfies that $\Theta'' H_{\theta_i}, 1 \leq i \leq 6$ belong to 
the ideal generated by $H_{\theta_j}, 1 \leq j \leq 6$, and that 
the integral curves of $\Theta''$ are singular curves with singular characteristics in 
$\{ H_{\theta_1} = \cdots = H_{\theta_6} = 0, H_{\theta_7} \not= 0\}$. 
Therefore the directions of singular curves of this case form 
the subbundle $\langle \theta_1, \theta_2\rangle$ generated by $\theta_1, \theta_2$. 
Therefore $\lambda_1$-$\lambda_2$-directions are singular directions in this case. 

Suppose $H_{\theta_6} = 0, H_{\theta_7} = 0$. Then $\lambda_3 H_{\theta_8} = 0, \lambda_1 H_{\theta_8} = 0$. 
Suppose $H_{\theta_8} \not= 0$. 
Then we have $H_{\theta_8} \not= 0$, and that $\lambda_1 = 0, \lambda_3 = 0$. 
We set $\Theta^{(3)} = \overrightarrow{H}_{\theta_2}$. Then we see that 
$\Theta''' H_{\theta_i}, 1 \leq i \leq 7$ belong to the ideal generated by $H_{\theta_j}, 1 \leq j \leq 7$, 
and that the integral curves of $\Theta^{(3)}$ are singular curves with singular characteristics in 
$F^{(3)} \setminus F^{(4)}$. In particular $\lambda_2$-axis is the singular direction in this case. 

Suppose $H_{\theta_8} = 0$. Then $\lambda_2 H_{\theta_9} = 0$. 
If $H_{\theta_9} = 0$, then the singular characteristic has to be in the zero-section of $T^*W$. 
Therefore, for a singular characteristic curve, we have to have $H_{\theta_8} \not= 0$, 
and then $\lambda_2 = 0$. 
Then again 
$\Theta^{(4)} = A\overrightarrow{H}_{\theta_1} + B\overrightarrow{H}_{\theta_3}$, where $A, B$ are arbitrary functions 
on $T^*W$ with $(A, B) \not= (0, 0)$, satisfies that $\Theta^{(4)} H_{\theta_i}, 1 \leq i \leq 8$ belong to 
the ideal generated by $H_{\theta_j}, 1 \leq j \leq 8$, and that 
the integral curves of $\Theta^{(4)}$ are singular curves with singular characteristics in 
$F^{(4)} \setminus F^{(5)}$. 
Therefore $\lambda_1$-$\lambda_3$-directions are singular directions in this case. 
\QED

Then we have 

\bep
\label{B_3(1,2,3)-and-(2,3)}
Let $F \subset TW$ be a germ of $(3,5,7,8,9)$-distribution with a pseudo-product 
structure $F = F_1 \oplus F_2 \oplus F_3$ of type $B_3(1,2,3)$ on a $9$-dimensional manifold $W$. 
Then the reduction $(Z', E')$ of $(W, F)$ by $F_1$ is a 
$(3,5,7,8)$-distribution with a pseudo-product structure $E' = E'_1\oplus E'_2$ of type $B_3(2,3)$. 
Moreover the prolongation $(W', F')$ of the reduction $(Z', E')$ by $E_2'$ is 
isomorphic, as a distribution with pseudo-product structure of type $B_3(1,2,3)$,  
to the original distribution $(W, F)$. 

Conversely, let $(Z, E)$ is a $(3,5,7,8)$-distribution with a pseudo-product structure $E = E_1\oplus E_2$ 
of type $B_3(2,3)$ and $(W, F)$ is the prolongation of $(Z, E)$ by $E_2$. 
Then the $(3,5,7,8,9)$-distribution $(W, F)$ with a pseudo-product 
structure $F = F_1 \oplus F_2 \oplus F_3$ of type $B_3(1,2,3)$
is reduced by $F_1$, to a $(3,5,7,8)$-distribution $(Z', E')$ with 
a pseudo-product structure $E' = E'_1\oplus E'_2$ is isomorphic to the original $(Z, E)$ with $E = E_1\oplus E_2$. 
\enp

\Proof
Let $F_i = \langle \theta_i\rangle, i = 1, 2, 3, F = \langle\theta_1, \theta_2, \theta_3\rangle$ satisfy the conditions 
of Definition \ref{B_3-(1,2,3)}. 
Then the distribution 
$\widetilde{E} = \langle \theta_1, \theta_2, \theta_3, \theta_4\rangle, 
\theta_4 := [\theta_1, \theta_2]$ is reduced by $F_1$ to 
the distribution $E' = \overline{\langle \theta_1, \theta_2, \theta_3, \theta_4\rangle} = \overline{\langle \theta_2, \theta_3, \theta_4\rangle}$
of rank $3$ on the leaf-space $Z' = W/F_1$, which is of dimension $8$. 
We have 
\[
F \subsetneq \widetilde{E} \subsetneq F^{(2)} \subsetneq \widetilde{E}^{(2)} \subsetneq F^{(3)} \subseteq \widetilde{E}^{(3)} = F^{(4)} \subsetneq \widetilde{E}^{(4)} = F^{(5)} = TW, 
\]
where $F^{(2)} = \langle \theta_i \mid 1 \leq i \leq 5\rangle, 
\widetilde{E}^{(2)} = \langle \theta_i \mid 1 \leq i \leq 6\rangle,
F^{(3)} = \langle \theta_i \mid 1 \leq i \leq 7\rangle, 
\widetilde{E}^{(3)} = \langle \theta_i \mid 1 \leq i \leq 8\rangle$, and 
$\widetilde{E}^{(4)} = \langle \theta_i \mid 1 \leq i \leq 9\rangle = TW$. 
Therefore we see that $E'$ is a $(3, 5, 7, 8)$-distribution. We set 
$E_1' := \overline{\langle \theta_1, \theta_3\rangle} = \overline{\langle \theta_3\rangle}, 
E_2' := \overline{\langle \theta_1, \theta_2, \theta_4\rangle} = \overline{\langle \theta_2, \theta_4\rangle}$. 
Then $E' =  E_1' \oplus E_2'$ satisfies
\[
[E_1', E_2'] = {E'}^{(2)}, \ [E_1', {E'}^{(2)}] = {E'}^{(3)}, \ [E_2', {E'}^{(2)}] = {E'}^{(2)}, \ [E_1', {E'}^{(3)}] = {E'}^{(3)}, \ [E_2', {E'}^{(3)}] = TZ'. 
\]
Thus $(Z', E'), E' = E_1' \oplus E_2'$ is a $(3, 5, 7, 8)$-distribution with pseudo-product structure of type 
$B_3(2,3)$. 
Moreover, regarding Proposition \ref{(3,5,7,8)-equivalence}, we take any generator $\zeta_1$ of $E_1'$ and 
any generators $\zeta_2, \zeta_3$ of $E_2'$ and we consider the prolongation $F'$ of $E'$ by $E_2$, which is 
generated by 
\[
\Theta_1 := \frac{\pa}{\pa w}, \quad \Theta_2 := \zeta_2 + w\zeta_3, \quad  
\Theta_3 := \zeta_1 + a\zeta_2 + b\zeta_3 + c\frac{\pa}{\pa w}. 
\]
Then $(W', F')$ is naturally identified with the original $(W, F)$. 

Then converse is clear. 
\QED

\section{Another prolongation of a $(3, 6)$-distribution.}
\label{Another-prolongation-of-(3, 6)-distribution}

Let $(X, D)$ be a $(3, 6)$-distribution and $D^* = {\mbox{\rm{Hom}}}(D, \R)$, the dual bundle of $D$. 
We set $S := P(D^*) = \{ (x, [\rho]) \mid x \in X, [\rho] \in P(D^*_x)\}$, which is identified with 
the space of $2$-dimensional subspaces in $D$ over any point of $X$. Note that $\dim(S) = 8$. 
Define the distribution $L \subset TS$ over $S$ by, for $(x, [\rho]) \in S$, 
\[
L_{(x, [\rho])} := \{ v \in T_{(x, [\rho])}S \mid \langle\rho, \pi_*(v)\rangle = 0 \}, 
\]
where $\pi : S \to X$ is the natural projection. Then $L$ is of rank $4$. 
Then we have 

\bep
\label{B_3(1,3)} 
The small growth vector of $L$ is given by $(4, 6, 8)$. Moreover 
$L$ has the decomposition $L = L_1\oplus L_2$ such that 
$L_1^{(2)}\oplus L_2$ is a pseudo-product structure of type $B_3(1,3)$ 
in the generalised sense of the following Definition \ref{pseudo-product-B_3(1,3)}. 
Furthermore the reduction of $(S, L)$ by $L_2$ is isomorphic to the original $(3, 6)$-distribution $(X, D)$. 
\enp

\bef
{\rm
\label{pseudo-product-B_3(1,3)}
Let $L \subset TS$ be a $(4, 6, 8)$-distribution with the decomposition $L = L_1\oplus L_2$ 
into subbundles of rank $2$. 
Then $L_1^{(2)}\oplus L_2$ is called 
a {\it pseudo-product structure of type $B_3(1,3)$ in the generalised sense} if 
there exist generators $\ell_1, \ell_2$ of $L_1$ 
and generators $\ell_3, \ell_4$ of 
$L_2$ respectively, satisfying that
\vspace{0.1truecm}
\\
$
{[\ell_1, \ell_2]} \equiv: \ell_5\ {\mbox{\rm mod.}}L_2, \ 
[\ell_1, \ell_3] \equiv: \ell_6\ {\mbox{\rm mod.}}L_2, \ [\ell_1, \ell_4] \in L_2, \ 
[\ell_2, \ell_3] \in L_2, \ [\ell_2, \ell_4] \equiv: \ell_6\ {\mbox{\rm mod.}}L_2, 
$
\\
$
{[\ell_3, \ell_4]} \in L_2, 
{\mbox{\rm so that\ }} L^{(2)} = \langle \ell_i \mid 1 \leq i \leq 6\rangle$, and $\rank(L^{(2)}) = 6, $
\\ 
${[\ell_1, \ell_5]} \in L^{(2)}, \ [\ell_1, \ell_6] \equiv: \ell_7, {\mbox{\rm mod.}}L^{(2)}, \ 
[\ell_2, \ell_5] \in L^{(2)}, \ [\ell_2, \ell_6] \equiv: \ell_8, {\mbox{\rm mod.}}L^{(2)}, 
[\ell_3, \ell_5] \equiv \ell_8,{\mbox{\rm mod.}}L^{(2)}, 
\\
{[\ell_3, \ell_6]} \in \langle \ell_1, \ell_2, \ell_3, \ell_4, \ell_6\rangle, \ 
[\ell_4, \ell_5] \equiv - \ell_7, {\mbox{\rm mod.}}L^{(2)}, \ 
[\ell_4, \ell_6] \in \langle \ell_1, \ell_2, \ell_3, \ell_4, \ell_6\rangle$, so that $L^{(3)} = \langle \ell_i \mid 1 \leq i \leq 8\rangle$ and $\rank(L^{(3)}) = 8. $
}
\enf

\ber
{\rm
Note that, by the reduction of $(S, L)$ to $(X, D)$ by $L_2$, 
the distribution $L_1$ projects to the $3$-dimensional cone field $\lambda(\xi_1 + y_1\xi_3) + \mu(\xi_2 + y_2\xi_3), 
(\lambda, \mu) \in \R^2, (y_1, y_2) \in (\R^2, 0)$, which has $D$ as the linear hull. 
}
\enr

\ber
\label{reason-to-add-generalised}
{\rm 
If $L_1^{(2)}\oplus L_2$ is a pseudo-product structure of type $B_3(1,3)$ in the generalised sense, 
then $L_2$ is integrable. However we do not suppose the integrability of $L_1^{(2)}$ 
in Definition \ref{pseudo-product-B_3(1,3)}. 
If also $L_1^{(2)}$ is integrable, then we call $L_1^{(2)}\oplus L_2$ 
a {\it pseudo-product structure of type $B_3(1,3)$ in the strict sense}, or simply 
a {\it pseudo-product structure of type $B_3(1,3)$}. 
Then there are induced locally double fibrations from $S$ to the $5$-dimensional leaf-space $S/L_1^{(2)}$ 
and to the $6$-dimensional leaf-space $S/L_2$. 
See \S \ref{B3-models} (iii), as an example of pseudo-product structures of type $B_3(1,3)$ in the strict sense, 
so that $L_1^{(2)}$ is integrable. 
Note that, even in the model case, $L_1$ is never integrable. This is the reason to call $L_1^{(2)}\oplus L_2$, 
not $L_1\oplus L_2$, 
a (generalised) pseudo-product structures of type $B_3(1,3)$. 
}
\enr

\noindent
{\it 
Proof of Proposition \ref{B_3(1,3)}.} \ 
We set $L_2 := \Ker(\pi_* : TS \to TX)$ globally 
and are going to fix the complementary bundle $L_1$ to $L_2$ in $L$. 

Take any $(x_0, [\rho_0]) \in S$. Without loss of generality, we assume $\Ker(\rho_0)$ projects 
isomorphically to $\langle \xi_1(x_0), \xi_2(x_0)\rangle \subset T_{x_0}X$. 
Take an open neighbourhood $U$ of $(x_0, [\rho_0])$ in $S$ such that, for any 
$(x, [\rho]) \in U$, $\Ker(\rho)$ projects to $\langle\xi_1(x), \xi_2(x)\rangle \subset T_xX$ isomorphically. 
Then we can take frame of $L\vert_U$ as 
\[
\ell_1 = \xi_1 + y_1\xi_3 + \alpha\frac{\pa}{\pa y_1} + \beta\frac{\pa}{\pa y_2}, \ \ 
\ell_2 = \xi_2 + y_2\xi_3 + \gamma\frac{\pa}{\pa y_1} + \delta\frac{\pa}{\pa y_2}, \ \ 
\ell_3 = \frac{\pa}{\pa y_1}, \ \ \ell_4 =  \frac{\pa}{\pa y_2}. 
\]
Here $\alpha, \beta, \gamma, \delta$ are some functions on $U$ which we fix later, 
and $y_1, y_2$ are inhomogeneous coordinates of $\pi$-fibres. Note that $L_2 = \langle \ell_3, \ell_4\rangle$. 

Then we have $[\ell_1, \ell_2] \equiv \xi_4 + y_2\xi_5 - y_1\xi_6 + (\beta - \gamma)\xi_3 =: \ell_5$, $[\ell_1, \ell_3] \equiv  - \xi_3 =: \ell_6$, 
$[\ell_1, \ell_4] \equiv 0, \ [\ell_2, \ell_3] \equiv 0, \ [\ell_2, \ell_4] \equiv - \xi_3 = \ell_6, \ 
[\ell_3, \ell_4] \equiv 0$, mod.$L_2$, therefore we have 
$L^{(2)} = \langle \ell_1, \ell_2, \ell_3, \ell_4, \ell_5, \ell_6\rangle$ and $\rank(L^{(2)}) = 6$. 

Moreover we have, modulo $L^{(2)}$, 
\[
\begin{array}{c}
{[\ell_1, \ell_5]} \equiv [\xi_1 + y_1\xi_3, \xi_4 + y_2\xi_5 - y_1\xi_6] + (2\beta - \gamma)\xi_5 - \alpha\xi_6, 
\\
{[\ell_2, \ell_5]} \equiv 
[\xi_2 + y_2\xi_3, \xi_4 + y_2\xi_5 - y_1\xi_6] + \delta\xi_5 + (\beta - 2 \gamma)\xi_6. 
\end{array}
\]
Since the quotient bundle $(TS/L^{(2)})\vert_U$ is generated by $\xi_5, \xi_6$, 
we see that there exist $(\alpha, \beta, \gamma, \delta)$ uniquely such that 
$[\ell_1, \ell_5] \equiv 0, [\ell_2, \ell_5] \equiv 0$, mod.$L^{(2)}$. 
Further we have
$[\ell_1, \ell_6] \equiv - \xi_5 =: \ell_7, \ 
[\ell_2, \ell_6] \equiv - \xi_6 =: \ell_8, \ 
[\ell_3, \ell_5] \equiv - \xi_6 = \ell_8, \ 
[\ell_4, \ell_5] \equiv \xi_5 = -\ell_7$, mod.$L^{(2)}$. 
Moreover we have both $[\ell_3, \ell_6], [\ell_4, \ell_6]$ belong to $L_2$. 
Therefore we have 
$L^{(3)} = \langle \ell_1, \ell_2, \ell_3, \ell_4, \ell_5, \ell_6, \ell_7, \ell_8\rangle = TU$ and we see that $\rank(L^{(3)}) = 8$. 

Now we set $L_1 := \langle \ell_1, \ell_2\rangle$. 
Then we have that $L$ has $(4, 6, 8)$ as the small growth vector and moreover the 
decomposition $L = L_1\oplus L_2$ gives 
a pseudo-product structure of type $B_3(1,3)$ in the generalised sense defined in Definition \ref{pseudo-product-B_3(1,3)}. 

For other coordinate neighbourhoods, the calculations are performed similarly. Note that 
the local uniqueness of $L_1$ implies the global existence and uniqueness of $L_1$. 

Lastly consider the local reduction of $(S, L)$ by the integrable distribution $L_2$. 
Return to the local calculations on $U$. 
Then we have $[L_1, L_2] = \langle \ell_1, \ell_2, \ell_3, \ell_4, \ell_6\rangle 
= \pi_*^{-1}D$, which reduces to $D$ by the reduction by $L_2$. 
\QED

\

Moreover conversely we have 

\bep
\label{B_3(1,3)-converse}
Starting from a $(4, 6, 8)$-distribution $(S, L)$, $L = L_1\oplus L_2$, with generalised pseudo-product structure 
$L_1^{(2)}\oplus L_2$ of type $B_3(1, 3)$, 
the local reduction of $(S, L)$ by $L_2$ turns to be a $(3, 6)$-distribution $(X', D')$. 
Moreover the prolongation $(S', L')$ of $(X', D')$ is isomorphic to $(S, L)$ 
as $(4, 6, 8)$-distributions with pseudo-product structure of type $B_3(1, 3)$. 
\enp

\Proof
Let $L_1 = \langle \ell_1, \ell_2\rangle, L_2 = \langle \ell_3, \ell_4\rangle$ and suppose 
the generators satisfy the conditions of Definition \ref{pseudo-product-B_3(1,3)}. 
Then $[L_1, L_2] = \langle \ell_1, \ell_2, \ell_3, \ell_4, \ell_6\rangle =: {\mathcal L}$. 
From $[\ell_3, \ell_6] \in {\mathcal L}, [\ell_4, \ell_6] \in {\mathcal L}$, we reduce $(S, L)$ 
as a local quotient, to $X' := S/L_2$ 
and $D' := {\mathcal L}/L_2 = \langle \overline{\ell}_1, \overline{\ell}_2, \overline{\ell}_6\rangle \subset TX'$, 
which is a distribution of rank $3$ on the $6$-dimensional space. 
Since $[\ell_1, \ell_2] \equiv \ell_5, [\ell_1, \ell_6] \equiv \ell_7, [\ell_2, \ell_6] \equiv \ell_8$ mod.$L_2$, 
we have 
$[{\mathcal L}, {\mathcal L}] = TS$ and therefore $[D', D'] = TX'$. Thus $(X', D')$ is a $(3, 6)$-distribution. 

We take a system of local coordinates $(x, y) = (x_1, \dots, x_6, y_1, y_2)$ of $S$ 
such that $\ell_3 = \frac{\pa}{\pa y_1}, \ell_4 = \frac{\pa}{\pa y_2}$. 
By the conditions that both $[\ell_3, \ell_6]$ and $[\ell_4, \ell_6]$ belong to 
$L_2 = \langle \ell_3, \ell_4\rangle$, we have 
$
\ell_1 \equiv \xi_1 + a\xi_3, 
\ 
\ell_2 \equiv \xi_2 + b\xi_3$, mod.$L_2 = \langle \ell_3, \ell_4\rangle$, where we set 
$\xi_1 = \overline{\ell}_1, \xi_2 = \overline{\ell}_2, \xi_3 = -\overline{\ell}_6$, and 
$a, b$ are functions on $(x, y)$. 
Since $[\ell_1, \ell_3] \equiv \frac{\pa a}{\pa y_1}\ell_6, [\ell_1, \ell_3] \equiv \ell_6$ mod.$L$, 
we have $a = y_1$. Similarly, since 
$[\ell_2, \ell_4] \equiv \frac{\pa b}{\pa y_2}\ell_6, [\ell_2, \ell_4] \equiv \ell_6$ mod.$L$, 
we have $b = y_2$. 
Thus we have $\ell_1 = \xi_1 + y_1\xi_3, \ell_2 = \xi_2 + y_1\xi_3$ with 
$\ell_3 = \frac{\pa}{\pa y_1}, \ell_4 = \frac{\pa}{\pa y_2}$. Hence we see that the prolongation $(S', L')$ 
of $(X, D)$ is isomorphic to $(S, L)$. 
\QED

\bee
\label{example-3} 
{\rm 
Returning to a $(3, 6)$-distribution $(X, D)$ in Example \ref{example-1}, 
we consider its prolongation $(S, L)$, $L = L_1\oplus L_2$ on $P(D^*)$ as in Proposition \ref{B_3(1,3)}.  
Then $(S, L)$ is given by $L = L_1\oplus L_2$, $L_1 = \langle \ell_1, \ell_2\rangle, L_2 = \langle \ell_3, \ell_4\rangle$, 
$\ell_1 = \xi_1 + y_1\xi_3 + \alpha\frac{\pa}{\pa y_1} + \beta\frac{\pa}{\pa y_2}, 
\ell_2 = \xi_2 + y_2\xi_3 + \gamma\frac{\pa}{\pa y_1} + \delta\frac{\pa}{\pa y_2}, \
\ell_3 = \frac{\pa}{\pa y_1}, \ell_4 =  \frac{\pa}{\pa y_2}$ with  
$\alpha = - y_1m'(x_6), \beta = - \frac{2}{3}y_2m'(x_6), \gamma = - \frac{1}{3}y_2m'(x_6), \delta = 0$, 
in this case. Then we see that $L_1^{(2)} = \langle \ell_1, \ell_2, [\ell_1, \ell_2]\rangle$ is never integrable provided 
$m(x_6)$ is not locally constant, 
so that $L_1^{(2)}\oplus L_2$ is pseudo-product structure of type $B_3(1,3)$ in the generalised sense 
but not in the strict sense. 
\qed
}
\ene

\section{Another iterated prolongation of a $(3, 6)$-distribution.} 
\label{another-iterated-prolongation}

Let us study singular curves of $(4, 6, 8)$-distributions $L = L_1\oplus L_2 \subset TS$ 
with generalized pseudo-product structure $L_1^{(2)}\oplus L_2$ 
introduced in \S \ref{Another-prolongation-of-(3, 6)-distribution} and 
try to prolong $(S, L)$ to get $(3, 5, 7, 8, 9)$-distributions. 
See Definition \ref{singular-velocity-cone} for singular velocity cone (SVC). 

\bel
\label{SVC(L)}
We have 
$
\SVC(L) = \SVC(L, L^{(2)\perp}) = 
\{ \sum_{i=1}^4 \lambda_i\ell_i \mid \lambda_1\lambda_3 + \lambda_2\lambda_4 = 0\}. 
$
Moreover we have $\SVC(L, L^\perp \setminus L^{(2)\perp}) = L_2$. 

\enl

\Proof\ 
The Hamiltonian of the control system defined by $(S, L)$ is given by 
\[
H = \lambda_1 H_{\ell_1} + \lambda_2 H_{\ell_2} + \lambda_3 H_{\ell_3} + \lambda_4 H_{\ell_4}, 
\]
where $H_{\ell_i}(x, y; p, q) = \langle (p, q), \ell_i(x, y)\rangle$ 
and $\lambda_1, \lambda_2, \lambda_3, \lambda_4$ are control parameters. 
The constrained Hamiltonian equation for singular curves of $(S, L)$ is given by 
\[
\dot{x} = \frac{\pa H}{\pa p}, \dot{y} = \frac{\pa H}{\pa q}, \dot{p} = - \frac{\pa H}{\pa x}, \dot{q} = - \frac{\pa H}{\pa y}, 
\]
with constraints $H_{\ell_1} = 0, H_{\ell_2} = 0, H_{\ell_3} = 0, H_{\ell_4} = 0$ and $(p, q) \not= (0, 0)$. 
Then along singular curves we have 
\[
0 = \frac{d}{dt} H_{\ell_1} = \lambda_1 H_{[\ell_1, \ell_1]} + \lambda_2 H_{[\ell_2, \ell_1]} + \lambda_3 H_{[\ell_3, \ell_1]} + \lambda_4 H_{[\ell_4, \ell_1]} = - \lambda_2 H_{\ell_5} - \lambda_3 H_{\ell_6}. 
\]
Similarly we have 
\[
\begin{array}{c}
0 = \frac{d}{dt} H_{\ell_2}  =  \lambda_1 H_{\ell_5} - \lambda_4 H_{\ell_6}, 
\quad
0 = \frac{d}{dt} H_{\ell_3}  =  \lambda_1 H_{\ell_6}, 
\quad
0 = \frac{d}{dt} H_{\ell_4}  =  \lambda_2 H_{\ell_6}, 
\end{array}
\]

Suppose $H_{\ell_6} \not= 0$. Then $\lambda_1 = 0, \lambda_2 = 0$, and then $\lambda_3 = 0, \lambda_4 = 0$, which 
implies the constancy of the singular curve. 

Let $H_{\ell_6} = 0$ and $H_{\ell_5} \not= 0$. Then $\lambda_1 = 0, \lambda_2 = 0$. 
Consider a vector field ${\mathcal V} = A\overrightarrow{H}_{\ell_3} + B\overrightarrow{H}_{\ell_4}$ 
over $T^*(U)$, for any function $A, B$ on $x, y, p, q$. 
Then we see all of functions ${\mathcal V}H_{\ell_i}, i = 1, 2, 3, 4, 6$ in $T^*U$ belong to ideal generated by 
functions $H_{\ell_j}, j = 1, 2, 3, 4, 6$ and the integral curves of ${\mathcal V}$ 
are abnormal bi-extremal (singular characteristic) 
of $L$. 
Then the corresponding singular velocities form $\SVC(L, L^\perp \setminus L^{(2)\perp}) = \langle \ell_3, \ell_4\rangle = L_2$. 

Now suppose that $H_{\ell_5} = 0, H_{\ell_6} = 0$, i.e., we treat the singular characteristic curves belonging to $L^{(2)\perp}$. 
Since
\[
\begin{array}{c}
\frac{d}{dt} H_{\ell_5}  =  \lambda_3 H_{\ell_8} - \lambda_4 H_{\ell_7}, 
\quad
\frac{d}{dt} H_{\ell_6}  =  \lambda_1 H_{\ell_7} + \lambda_2 H_{\ell_8}, 
\end{array}
\]
we have 
\[
\left(
\begin{array}{cc}
-\lambda_4 & \lambda_3
\\
\lambda_1 & \lambda_2
\end{array}
\right)
\left(
\begin{array}{c}
H_{\lambda_7}
\\
H_{\lambda_8}
\end{array}
\right)
= 
\left(
\begin{array}{c}
0
\\
0
\end{array}
\right). 
\]
From $(H_{\lambda_7}, H_{\lambda_8}) \not= (0, 0)$, we have $\lambda_1\lambda_3 + \lambda_2\lambda_4 = 0$. 
Regarding these necessary conditions on singular curves, 
we set
\[
\widetilde{\mathcal{V}} = \widetilde{A}(-H_{\ell_8}\overrightarrow{H}_{\ell_1} + H_{\ell_7}\overrightarrow{H}_{\ell_2}) 
+ \widetilde{B}(H_{\ell_7}\overrightarrow{H}_{\ell_3} + H_{\ell_8}\overrightarrow{H}_{\ell_4}), 
\]
where $\widetilde{A}$ and $\widetilde{B}$ are $C^\infty$ functions on $T^*(U)$. 
Then we see that, for $1 \leq i \leq 6$, $\widetilde{\mathcal{V}}H_{\ell_i}$ belongs to 
the ideal generated by $H_{\ell_j}, 1 \leq j \leq 6$ over the algebra of $C^\infty$ functions on $T^*(U)$. 
Therefore we have that a curve is a singular characteristic of $L$ if and only if it is an integral curve of a vector field 
of type ${\mathcal V}$ or $\widetilde{\mathcal V}$. 
In particular we have that 
$\SVC(L, L^{(2)\perp}) = 
\{ \sum_{i=1}^4 \lambda_i\ell_i \mid \lambda_1\lambda_3 + \lambda_2\lambda_4 = 0\}$. 
Thus we have the result. 
\QED

\

Now we prolong $(S, L)$ regarding the singularities of $L = L_1\oplus L_2 = 
\langle \ell_1, \ell_2\rangle\oplus\langle\ell_3, \ell_4\rangle$. 
Recall that $\ell_1 = \xi_1 + y_1\xi_3 + \alpha\frac{\pa}{\pa y_1} + \beta\frac{\pa}{\pa y_2}, 
\ell_2 =  \xi_2 + y_2\xi_3 + \gamma\frac{\pa}{\pa y_1} + \delta\frac{\pa}{\pa y_2}$, and 
$\alpha, \beta, \gamma, \delta$ are fixed already in the proof of Lemma \ref{B_3(1,3)}. 

We take any subbundle $C \subset SVC(L)$ of rank $2$ which satisfies 
$C \cap L_1 = \{ 0\}$ and $C \cap L_2 = \{ 0\}$. 
By changing $\ell_1, \ell_2, \ell_3, \ell_4$ by their functional multiples respectively, 
we put 
\[
C := \{ \lambda_1\ell_1(s) + \lambda_2\ell_2(s) - \lambda_2\ell_3(s) + \lambda_1\ell_4(s) \mid 
s \in S, \lambda_1, \lambda_2 \in \R \},
\]
without loss of generalities in the following discussion. 
We set $\widetilde{W} := P(C)$, the space of linear lines in $C$. 
Denote by $\pi_S : \widetilde{W} \to S$ the canonical projection and define $\widetilde{F} \subset T\widetilde{W}$ by 
\[
\widetilde{F}_{(s, [\lambda])} := \{ v \in T_{(s, [\lambda])}\widetilde{W} \mid \pi_{S*}(v) 
\in \lambda_1\ell_1(s) + \lambda_2\ell_2(s) - \lambda_2\ell_3(s) + \lambda_1\ell_4(s)  \}, 
\]
for $(s, [\lambda]) \in \widetilde{W}$. Then we see $\dim(\widetilde{W}) = 9$. 

In the coordinate neighbourhood $\lambda_1 \not= 0$ of $\widetilde{W}$, 
we set $u = \lambda_2/\lambda_1$. 
Then we have the local frame of $\widetilde{F}$ : 
\[
\tau_1 = \ell_4 - u\ell_3, \tau_2 = \frac{\pa}{\pa u}, \tau_3 = \ell_1 + u\ell_2 + \eta(\ell_4 - u\ell_3) 
+ \varepsilon\frac{\pa}{\pa u},  
\]
where $\{\ell_1, \ell_2, \ell_3, \ell_4\}$ is a local frame given in the proof of Proposition \ref{B_3(1,3)} and 
$\eta, \varepsilon$ are functions on $\widetilde{W}$ which will be fixed later. 
We take a generalised Legendre transformation $y_1 = z_3 - wz_2, \ y_2 = w, \ u = z_2$ on $\widetilde{W}$ from 
$\tau_1, \tau_2, \tau_3$ (see Remark \ref{generalised-Legendre-transformation}) to get 
\[
\begin{array}{rcl}
\tau_1 & = & \ell_4 - u\ell_3 \ = \ \frac{\pa}{\pa w}, 
\qquad 
\tau_2  = \frac{\pa}{\pa u} \ = \ \frac{\pa}{\pa z_2} + w\frac{\pa}{\pa z_3}, 
\vspace{0.2truecm}
\\
\tau_3 & = & \ell_1 + u\ell_2 + \eta(\ell_4 - u\ell_3) + \varepsilon\frac{\pa}{\pa u} 
\vspace{0.1truecm}
\\
& = & 
\xi_1 + (z_3 - wz_2)\xi_3 + \alpha\frac{\pa}{\pa z_3} + \beta(\frac{\pa}{\pa w} + z_2\frac{\pa}{\pa z_3}) 
\vspace{0.1truecm}
\\
& & 
\hspace{3truecm}
+ z_2\{ \xi_2 + w\xi_3 + \gamma\frac{\pa}{\pa z_3} + \delta(\frac{\pa}{\pa w} + z_2\frac{\pa}{\pa z_3})\} 
+ \eta\frac{\pa}{\pa u} 
+ \varepsilon(\frac{\pa}{\pa z_2} + w\frac{\pa}{\pa z_3})
\vspace{0.2truecm}
\\
& = & 
\xi_1 + z_2\xi_2 + z_3\xi_3 + 
A\frac{\pa}{\pa z_2} + B\frac{\pa}{\pa z_3} + C\frac{\pa}{\pa w}, 
\end{array}
\]
where $A = \varepsilon, B = \alpha + z_2(\beta + \gamma) + z_2^2\delta + w\varepsilon$ and 
$C = \eta + \beta + z_2\delta$. 
Moreover we set $\widetilde{F}_1 = \langle \tau_1\rangle, \widetilde{F}_2 = \langle \tau_2\rangle$ 
and $\widetilde{F}_3  = \langle \tau_3\rangle$.  
Note that $\varepsilon, \eta$ and so $A, B, C$ are uniquely determined in the proof of the following: 

\bep
\label{Another-B_3(1,2,3)}
The small growth vector of $\widetilde{F}$ is given by $(3, 5, 7, 8, 9)$. 
Moreover 
$\widetilde{F} = \widetilde{F}_1\oplus\widetilde{F}_2\oplus\widetilde{F}_3$ is a pseudo-product structure 
of type $B_3(1,2,3)$ in the sense of Definition \ref{B_3-(1,2,3)}, if and only if $A$ satisfies that 
$A_w = 0$ and 
$[\xi_1 + z_2\xi_2 + z_3\xi_3, \xi_5 + z_2\xi_6] + A\xi_6 \in \widetilde{F}^{(4)}$. 
\enp

\Proof
We have 
$[\tau_1, \tau_2] = \frac{\pa}{\pa z_3} =: \tau_4, [\tau_1, \tau_3] =  
A_w\frac{\pa}{\pa z_2} + B_w\frac{\pa}{\pa z_3} + C_w\frac{\pa}{\pa w} \equiv (B_w - wA_w)\frac{\pa}{\pa z_3}$ 
mod.$\widetilde{F}$. Since $B_w - wA_w = \{\alpha + z_2(\beta + \gamma) + z_2^2\delta\}_w+ \varepsilon$, 
we see $[\tau_1, \tau_3] \equiv 0$ mod.$\widetilde{F}$, 
if we take $\varepsilon = A = -\{\alpha + z_2(\beta + \gamma) + z_2^2\delta\}_w$. 

Remark that $[\tau_1, \tau_3] \in \langle \tau_1, \tau_3\rangle$ if and only if 
$A_w = 0, B_w = 0$. Note that $A_w = 0$ implies $B_w = 0$. 

Further we have $[\tau_2, \tau_3] \equiv 
\xi_2 + w\xi_3 + \tau_2(A)\frac{\pa}{\pa z_2} + \{\tau_2(B) - C\}\frac{\pa}{\pa z_3} + \tau_3(C)\frac{\pa}{\pa w} 
=: \tau_5$, mod.$\widetilde{F}$, and 
$\widetilde{F}^{(2)} = \langle \tau_i \mid 1 \leq i \leq 5\rangle$, $\rank(\widetilde{F}^{(2)}) = 5$. 

Moreover we have 
$[\tau_1, \tau_4] = 0, [\tau_1, \tau_5] \equiv \xi_3 =: \tau_6$, 
$[\tau_2, \tau_4] = 0, [\tau_2, \tau_5] \equiv 0$, 
$[\tau_3, \tau_4] \equiv - \xi_3 = - \tau_6, 
[\tau_3, \tau_5] \equiv 
\xi_4 - z_3\xi_6 + w(\xi_5 + z_2\xi_6) - \tau_2(A)\xi_2 - \tau_2(B)\xi_3 + C\xi_3 =: \tau_7$, mod.$\widetilde{F}^{(2)}$, and 
$\widetilde{F}^{(3)} = \langle \tau_i \mid 1 \leq i \leq 7\rangle \supset \langle \frac{\pa}{\pa z_2}, \frac{\pa}{\pa z_3}, 
\xi_1, \xi_2, \xi_3\rangle$, $\rank(\widetilde{F}^{(3)}) = 7$. 

Furthermore we have 
$[\tau_1, \tau_6] \equiv 0, [\tau_1, \tau_7] \equiv \xi_5 + u\xi_6 =: \tau_8, 
[\tau_2, \tau_6] \equiv 0, [\tau_2, \tau_7] \equiv 0, 
[\tau_3, \tau_6] \equiv \xi_5 + z_2\xi_6 = \tau_8$, and 
\[
\begin{array}{rcl}
[\tau_3, \tau_7] & \equiv & 
[\xi_1 + z_2\xi_2 + z_3\xi_3, \xi_4 - z_3\xi_6] - B\xi_6 
+ w([\xi_1 + z_2\xi_2 + z_3\xi_3, \xi_5 + z_2\xi_6] + A\xi_6) 
\\
& & \qquad\qquad\qquad\qquad + (C - \tau_2(B))(\xi_5 + z_2\xi_6) 
- \tau_2(A)(\xi_4 - z_3\xi_6) 
\\
& \equiv & 
{[\xi_1 + z_2\xi_2 + z_3\xi_3, \xi_4 - z_3\xi_6]} - B\xi_6 
+ w([\xi_1 + z_2\xi_2 + z_3\xi_3, \xi_5 + z_2\xi_6] + A\xi_6) 
\\
& \equiv & - \tau_2(A)\{\xi_4 - z_3\xi_6 + w(\xi_5 + z_2\xi_6)\} + 
\{ C - \tau_2(B) + w\tau_2(A)\}(\xi_5 + z_2\xi_6)
\\
& \equiv&  0\ \ {\mbox{\rm mod.}}\ \widetilde{F}^{(3)}, 
\end{array}
\]
if we set $\eta = \tau_2(B) - w\tau_2(A) - \beta - z_2\delta$. Thus we have 
$\widetilde{F}^{(4)} = \langle \tau_i \mid 1 \leq i \leq 8\rangle$ and $\rank(\widetilde{F}^{(4)}) = 8$. 

Moreover we have $[\tau_1, \tau_8] \equiv 0$ mod.$\widetilde{F}^{(4)}$, 
$[\tau_2, \tau_8] \equiv \xi_6 =: \tau_9$ mod.$\widetilde{F}^{(4)}$, and 
$[\tau_3, \tau_8] \equiv [\xi_1 + z_2\xi_2 + z_3\xi_3, \xi_5 + z_2\xi_6] + A\xi_6$ mod.$\widetilde{F}^{(4)}$.  
Thus we see that $\widetilde{F}^{(5)} = \langle \tau_i \mid 1 \leq i \leq 9\rangle = T\widetilde{W}$ 
and $\rank(\widetilde{F}^{(5)}) = 9$.  Note that 
$[\tau_3, \tau_8]$ needs not to be contained in $\widetilde{F}^{(4)}$. 

Thus we have the small growth vector of $\widetilde{F}$ is equal to $(3, 5, 7, 8, 9)$. 

Regarding Definition \ref{B_3-(1,2,3)}, we have the rest of the claim. 
\QED

\ber
\label{generalised-Legendre-transformation}
{\rm 
For the comparison of two classes of $(3, 5, 7, 8, 9)$-distributions $(W, F)$ and $(\widetilde{W}, \widetilde{F})$ given 
in Proposition \ref{iterated-prolongation} and Proposition \ref{Another-B_3(1,2,3)} respectively, 
we have used the transformation 
\[
y_1 = z_3 - wz_2, \ y_2 = w, \ u = z_2, \quad   {\mbox{\rm or equivalently}} \quad 
w = y_2, \ z_2 = u, \ z_3 = y_1 + uy_2, 
\] 
which remind us of the Legendre transformation in contact geometry (see \cite{Arnold} for instance), 
and can be called a {\it generalised Legendre transformation}. 
For the $B_3$-model (see \S \ref{B3-models} (i)), the transformation 
provides the same presentation of the pseudo-product structure of type $B_3(1,2,3)$ (see Example \ref{example-4} below). 
}
\enr

\bee
\label{example-4}
{\rm
Recall Examples \ref{example-1}, \ref{example-2} and \ref{example-3}. 
For those examples, we calculate functions $\alpha, \beta, \gamma, A, B, C$ appeared 
in the proof of Proposition \ref{Another-B_3(1,2,3)} as follows: 
\[
\begin{array}{c}
\alpha = - (z_3 - wz_2)\,m'(x_6), \ \  \beta = - \frac{2}{3}w\,m'(x_6), \ \ \gamma = - \frac{1}{3}w\,m'(x_6), \ \ \delta = 0, 
\vspace{0.2truecm}
\\
A = 0, \ \  B = - z_3\,m'(x_6), \ \ C = - w\,m'(x_6). 
\end{array}
\]
In particular, we have that $A_w = 0$. 
Moreover we see that $[\tau_3, \tau_8] \equiv 0$ mod.$\widetilde{F}^{(3)}$ and 
that another prolongation $(\widetilde{W}, \widetilde{F})$ 
coincides with $(W, F)$ in Proposition \ref{iterated-prolongation} 
in the case of Example \ref{example-1} in fact. 
}
\ene

\section{Proof of the main theorem.}
\label{Proof of the main theorem}

In this section we prove the main Theorem \ref{main-theorem}. 

By Proposition \ref{B_3(2,3)-and-(3,6)}, we see the local isomorphism classes of $(3, 6)$-distributions 
and 
those of $(3, 5, 7, 8)$-distributions with pseudo-product structure of type $B_3(2,3)$ are 
bijectively corresponds by prolongation and reduction procedures. 

By Proposition \ref{B_3(1,2,3)-and-(2,3)}, we see 
the local isomorphism classes of $(3, 5, 7, 8)$-distributions with pseudo-product structure of type $B_3(2,3)$ 
and 
those of $(3, 5, 7, 8, 9)$-distributions with pseudo-product structure of type $B_3(1,2,3)$ are 
bijectively corresponds by prolongation and reduction procedures. 

Further, by Propositions \ref{B_3(1,3)} and \ref{B_3(1,3)-converse}, 
the local isomorphism classes of $(3, 6)$-distributions corresponds bijectively 
to those of $(4, 6, 8)$-distributions with pseudo-product structure of type $B_3(1,3)$ in the generalised sense. 

These arguments prove Theorem \ref{main-theorem}. 
\QED

\section{Appendix: $B_3$-models.}
\label{B3-models}

Let $\R^{3,4}$ denote the vector space $\R^7 = \langle e_i \mid i = 1, \dots, 7\rangle$ with the 
inner product $(\ \mid \ )$ of index $(3,4)$ defined by 
\[
(e_i \vert e_j) = 
\left\{
\begin{array}{rl}
\frac{1}{2} & (i + j = 8, i \not= j), 
\vspace{0.1truecm}
\\
-\frac{1}{2} & (i = 4, j = 4), 
\\
0 & ({\mbox{\rm otherwise}}). 
\end{array}
\right.
\]

We consider spaces ${\mathscr F}_{1,2,3}, {\mathscr F}_{2,3}, {\mathscr F}_{1,3}, {\mathscr F}_3$ 
of flags by null subspaces in $\R^{3,4}$; 
\[
{\mathscr F}_{1,2,3} := \{ (V_1, V_2, V_3) \mid V_1 \subset V_2 \subset V_3 \subset \R^{3,4}, 
V_i {\mbox{\rm\ is a null subspace of dimension}}\ i, i = 1,2,3\}, 
\]
\[
{\mathscr F}_{2,3} := \{ (V_2, V_3) \mid V_2 \subset V_3 \subset \R^{3,4}, V_i {\mbox{\rm\ 
is a null subspace of dimension}}\ i, i = 2,3\}, 
\]
\[
{\mathscr F}_{1,3} := \{ (V_1, V_3) \mid V_1 \subset V_3 \subset \R^{3,4}, V_j {\mbox{\rm\ 
is a null subspace of dimension}}\ j, j = 1,3\}, 
\]
\[
{\mathscr F}_3 := \{ V_3 \mid V_3 \subset \R^{3,4}, {\mbox{\rm\ a null subspace of dimension}}\ 3\}. 
\]

The group ${\rm SO}(3, 4)$ of orientation-preserving linear isometries on $\R^{3,4}$ 
acts on each of the above spaces naturally and transitively. 
Denoted by $\pi_1 : {\mathscr F}_{1,2,3} \to {\mathscr F}_{2,3}$, $\pi_1 : {\mathscr F}_{1,3} \to {\mathscr F}_{3}$, 
$\pi_2 : {\mathscr F}_{1,2,3} \to {\mathscr F}_{1,3}$ and $\pi_2 : {\mathscr F}_{2,3} \to {\mathscr F}_{3}$ by 
the natural projections defined by forgetting $i$-th component for $i = 1, 2$, respectively. 

\

(i) (The $B_3(1,2,3)$-model) \ 
First we consider the space $W = {\mathscr F}_{1,2,3}$ consisting of complete null flags in $\R^{3,4}$. 
Let $(V_1^0, V_2^0, V_3^0) \in {\mathscr F}_{1,2,3}$, where $V_1^0 = \langle e_1\rangle, 
V_2^0 = \langle e_1, e_2\rangle, V_3^0 = \langle e_1, e_2, e_3\rangle$. 
Let $p :  \R^{3,4} \to V_3^0$ be the projection defined by 
$p(\sum_{i=1}^7 a_ie_i) := \sum_{i=1}^3a_ie_i$ and $U$ 
the set of flags  $(V_1, V_2, V_3) \in {\mathscr F}_{1,2,3}$ such that 
$p\vert_{V_3}$ is isomorphic. 
Take any flag 
$(V_1, V_2, V_3) \in U$ with $V_1 = \langle f_1\rangle, 
V_2 = \langle f_1, f_2\rangle, V_3 = \langle f_1, f_2, f_3\rangle$ and 
\[
\left\{ 
\begin{array}{rcrrrrrrl}
f_1 & = & 
e_1 & + w_{21}e_2 & + w_{31}e_3 & + w_{41}e_4 + w_{51}e_5 + w_{61}e_6 + w_{71}e_7, 
\\
f_2 & = & 
& e_2 & + w_{32}e_3 & + w_{42}e_4 + w_{52}e_5 + w_{62}e_6 + w_{72}e_7, 
\\
f_3 & = & 
 &     &  e_3 & + w_{43}e_4 + w_{53}e_5 + w_{63}e_6 + w_{73}e_7. 
\end{array}
\right.
\]
Then we have
\[
\begin{array}{rcl}
(f_1\mid f_1) & = & w_{71} + w_{21}w_{61} + w_{31}w_{51} - \frac{1}{2}w_{41}^2 = 0, 
\vspace{0.1truecm}
\\
(f_1\mid f_2) & = & \frac{1}{2}w_{72} + \frac{1}{2}w_{21}w_{62} + \frac{1}{2}w_{31}w_{52} - \frac{1}{2}w_{41}w_{42} + \frac{1}{2}w_{51}w_{32} + \frac{1}{2}w_{61} = 0, 
\vspace{0.1truecm}
\\
(f_1\mid f_3) & = & \frac{1}{2}w_{73} + \frac{1}{2}w_{21}w_{63} + \frac{1}{2}w_{31}w_{53} - \frac{1}{2}w_{41}w_{43} + \frac{1}{2}w_{51} = 0, 
\vspace{0.1truecm}
\\
(f_2\mid f_2) & = & w_{62} + w_{32}w_{52} - \frac{1}{2}w_{42}^2 = 0, 
\vspace{0.1truecm}
\\
(f_2\mid f_3) & = & \frac{1}{2}w_{63} + \frac{1}{2}w_{32}w_{53} - \frac{1}{2}w_{42}w_{43} + 
\frac{1}{2}w_{52} = 0, 
\vspace{0.1truecm}
\\
(f_3\mid f_3) & = & w_{53} - \frac{1}{2}w_{43}^2 = 0. 
\end{array}
\]
Then we have that $w_{21}, w_{31}, w_{41}, w_{51}, w_{61}, w_{32}, w_{42}, w_{52}, w_{43}$ 
form a system of coordinates on $U \subset {\mathscr F}_{1,2,3}$. 
The canonical differential system $F \subset T({\mathscr F}_{1,2.3})$ is defined by the infinitesimal condition 
\[
f_1'(t) \in V_2(t), \ f_2'(t) \in V_3(t), \ f_3'(t) \in V_3(t)^\perp \quad \cdots\cdots\cdots (\#), 
\]
at $t = 0$, for a curve $(V_1(t), V_2(t)), V_3(t)) = (\langle f_1(t)\rangle, \langle f_1(t), f_2(t)\rangle, \langle f_1(t), f_2(t), f_3(t)\rangle)$ on ${\mathscr F}_{1,2.3}$, where $V_3(t)^\perp$ means the orthogonal space to $V_3(t)$ for 
the indefinite metric $(\ \mid \ )$. 
From the conditions $f_1'(t) \in V_2(t), f_2'(t) \in V_3(t), (f_1\mid f_3) = 0, (f_2\mid f_3) = 0, (f_3\mid f_3) = 0$ at $t = 0$, 
we have $f_3' \in V_3(t)^\perp$ at $t = 0$, so the condition $(\#)$ is equivalent to, at $t = 0$, 
\[
f_1'(t) \in V_2(t), \ f_2'(t) \in V_3(t) \quad \cdots\cdots\cdots (\dag). 
\]
We see, by straightforward calculations, the condition $(\dag)$ is reduced to the differential system
\[
\begin{array}{c}
dw_{31} - w_{32}dw_{21} = 0, \ \ 
dw_{41} - w_{42}dw_{21} = 0, \ \ 
dw_{51} - w_{52}dw_{21} = 0, 
\\
dw_{61} + (w_{32}w_{52} - \frac{1}{2}w_{42}^2)dw_{21} = 0, \ \ 
dw_{42} - w_{43}dw_{32} = 0, \ \ 
dw_{52} - \frac{1}{4}w_{43}^2dw_{32} = 0, 
\end{array}
\]
which defines the distribution $F$ of rank $3$ generated by 
\[
\left\{ 
\begin{array}{rcl} 
\theta_1 & = & 
\frac{\pa}{\pa w_{21}} + w_{32}\frac{\pa}{\pa w_{31}} + w_{42}\frac{\pa}{\pa w_{41}} + w_{52}\frac{\pa}{\pa w_{51}} + \left( - w_{32}w_{52} + \frac{1}{2}w_{42}^2\right)\frac{\pa}{\pa w_{61}}, 
\vspace{0.2truecm}
\\
\theta_2 & = &  
\frac{\pa}{\pa w_{32}} + w_{43}\frac{\pa}{\pa w_{42}} + \frac{1}{2}w_{43}^2\frac{\pa}{\pa w_{52}}, 
\vspace{0.2truecm}
\\
\theta_3 & = & 
\frac{\pa}{\pa w_{43}}. 
\end{array}
\right.
\]
Then we have 
$[\theta_1, \theta_2] = - \frac{\pa}{\pa w_{31}} - w_{43}\frac{\pa}{\pa w_{41}} - \frac{1}{2}w_{43}^2\frac{\pa}{\pa w_{51}} 
+ (w_{52} - w_{42}w_{43} + \frac{1}{2}w_{32}w_{43}^2)\frac{\pa}{\pa w_{61}} =: \theta_4, 
[\theta_1, \theta_3] = 0, [\theta_2, \theta_3] = - \frac{\pa}{\pa w_{42}} - w_{43}\frac{\pa}{\pa w_{52}} =: \theta_5$, 
and $\rank(F^{(2)}) = 5$. Further we have 
$[\theta_1, \theta_4] = 0, 
[\theta_1, \theta_5] = \frac{\pa}{\pa w_{41}} + w_{43}\frac{\pa}{\pa w_{51}} + (w_{42} - w_{32}w_{43})\frac{\pa}{\pa w_{61}} =: \theta_6, 
[\theta_2, \theta_4] = 0,
[\theta_2, \theta_5] = 0, [\theta_3, \theta_4] = - \theta_6,  [\theta_3, \theta_5] = - \frac{\pa}{\pa w_{52}} =: \theta_7$, 
and thus $\rank(F^{(3)}) = 7$. Furthermore we have 
$[\theta_1, \theta_6] = 0, [\theta_1, \theta_7] = \frac{\pa}{\pa w_{51}} - w_{32}\frac{\pa}{\pa w_{61}} =: \theta_8, 
[\theta_2, \theta_6] = 0, [\theta_2, \theta_7] = 0, [\theta_3, \theta_6] = \theta_8, [\theta_3, \theta_7] = 0$, and 
$\rank(F^{(4)}) = 8$. 
Moreover $[\theta_1, \theta_8] = 0, [\theta_2, \theta_8] = - \frac{\pa}{\pa w_{61}} =: \theta_9, [\theta_3, \theta_8] = 0$, 
and therefore $\rank(F^{(5)}) = 9$ i.e., $F^{(5)} = T{\mathscr F}_{1,2,3}$. 
Thus we see $F$ has $(3, 5, 7, 8, 9)$ as the small growth vector. 

We set $F_1 = \langle\theta_1\rangle, F_2 = \langle\theta_2\rangle, F_3 = \langle\theta_3\rangle$. Then
$F = F_1\oplus F_2 \oplus F_3$ is a pseudo-product structure of type $B_3(1,2,3)$. 

The gradation Lia algebra is given, denoting the class of $\theta_i$ by $\delta_i$,  as 
\[
\begin{array}{c}
F \oplus F^{(2)}/F \oplus F^{(3)}/F^{(2)} \oplus F^{(4)}/F^{(3)} \oplus F^{(5)}/F^{(4)}, 
\\
F = \langle \delta_1, \delta_2, \delta_3\rangle, \ F^{(2)}/F = \langle \delta_4, \delta_5\rangle, \ 
F^{(3)}/F^{(2)} = \langle \delta_6, \delta_7\rangle, \ 
F^{(4)}/F^{(3)} = \langle \delta_8\rangle, \ F^{(5)}/F^{(4)} = \langle \delta_9\rangle, 
\end{array}
\]
with generators $\delta_1, \dots, \delta_9$ and operations 
$[\delta_1, \delta_2] = \delta_4, [\delta_1, \delta_3] = 0, [\delta_2, \delta_3] = \delta_5, 
[\delta_1, \delta_4] = 0, [\delta_1, \delta_5] = \delta_6, [\delta_2, \delta_4] = 0, 
[\delta_2, \delta_5] = 0, [\delta_3, \delta_4] = - \delta_6, [\delta_3, \delta_5] = \delta_7, 
[\delta_1, \delta_6] = 0, [\delta_1, \delta_7] = \delta_8, 
[\delta_2, \delta_6] = 0, [\delta_2, \delta_7] = 0, [\delta_3, \delta_6] = \delta_8, [\delta_3, \delta_7] = 0, 
[\delta_1, \delta_8] = 0, [\delta_2, \delta_8] = \delta_9, [\delta_3, \delta_8] = 0$. 

\

(ii) (The $B_3(2,3)$-model) \ 
Second we study $Z = {\mathscr F}_{2,3}$ which consists of null flags $V_2 \subset V_3 \subset \R^{3,4}$. 
Let $V_2 = \langle g_1, g_2\rangle, V_3 = \langle g_1, g_2, g_3\rangle$ with 
\[
\left\{ 
\begin{array}{rcrrrrrrl}
g_1 & = & 
e_1 &  & + z_{31}e_3 & + z_{41}e_4 + z_{51}e_5 + z_{61}e_6 + z_{71}e_7, 
\\
g_2 & = & 
& e_2 & + z_{32}e_3 & + z_{42}e_4 + z_{52}e_5 + z_{62}e_6 + z_{72}e_7, 
\\
g_3 & = & 
 &     &  e_3 & + z_{43}e_4 + z_{53}e_5 + z_{63}e_6 + z_{73}e_7. 
\end{array}
\right.
\]
The condition that $(V_2, V_3)$ is a null flag is given by
\[
\begin{array}{rcl}
(g_1 \vert g_1) & = & z_{71} + z_{31}z_{51} - \frac{1}{2}z_{41}^2 = 0, 
\\
(g_1 \vert g_2) & = & \frac{1}{2}z_{72} + \frac{1}{2}z_{31}z_{52} - \frac{1}{2}z_{41}z_{42} + \frac{1}{2}z_{32}z_{51} + \frac{1}{2}z_{61} = 0, 
\\
(g_1 \vert g_3) & = & \frac{1}{2}z_{73} + \frac{1}{2}z_{31}z_{53} - \frac{1}{2}z_{41}z_{43} + \frac{1}{2}z_{51} = 0, 
\\
(g_2 \vert g_2) & = & z_{62} + z_{32}z_{52} - \frac{1}{2}z_{42}^2 = 0, 
\\
(g_2 \vert g_3) & = & \frac{1}{2}z_{63} + \frac{1}{2}z_{32}z_{53} - \frac{1}{2}z_{42}z_{43} + \frac{1}{2}z_{52} = 0, 
\\
(g_3 \vert g_3) & = & z_{53} - \frac{1}{2}z_{43}^2 = 0. 
\end{array}
\]
In particular we see that $z_{31}, z_{41}, z_{51}, z_{61}, z_{32}, z_{42}, z_{52}$ and $z_{43}$ form 
a local system of coordinates of the null flag manifold ${\mathscr F}_{2,3}$. 
We define $E \subset T{\mathscr F}_{2,3}$ by the condition 
that $v \in E_{(V_2,V_3)} \subset T_{(V_2,V_3)}{\mathscr F}_{2,3}$ if and only if there is a curve $(V_2(t), V_3(t))$ 
in ${\mathscr F}_{2,3}$ satisfying $(V_2(0), V_3(0)) = (V_2, V_3)$, 
$v = (V_2'(0), V_3'(0))$, $V_2'(0) \in V_3$ and $V_3'(0) \in V_3^\perp$. 
We see the condition is equivalent, in terms of frames, $V_2(t) = \langle g_1(t), g_2(t)\rangle, V_3(t) = \langle 
g_1(t), g_2(t), g_3(t)\rangle$, to that 
$g_1'(0), g_2'(0) \in V_3$. Note that the condition $g_3'(0) \in V_3^\perp$ follows from that $g_1'(0), g_2'(0) \in V_3$. 
Then the condition is written by the Pfaff system
\[
\begin{array}{c}
dz_{41} - z_{43}dz_{31} = 0, \ dz_{51} - \frac{1}{2}z_{43}^2 dz_{31} = 0, \ dz_{61} + (\frac{1}{2}z_{32}z_{43}^2 - z_{42}z_{43} + z_{52})dz_{31} = 0, 
\vspace{0.2truecm}
\\
dz_{42} - z_{43}dz_{32} = 0, \ dz_{52} - \frac{1}{2} z_{43}^2 dz_{32} = 0, 
\end{array}
\] 
and we have
\[
\left\{ 
\begin{array}{rcl}
\zeta_1 & = & \frac{\pa}{\pa z_{43}}, 
\vspace{0.2truecm}
\\
\zeta_2 & = & \frac{\pa}{\pa z_{31}} + z_{43}\frac{\pa}{\pa z_{41}} + \frac{1}{2} z_{43}^2\frac{\pa}{\pa z_{51}}
- \left(\frac{1}{2}z_{32}z_{43}^2 - z_{42}z_{43} + z_{52}\right)\frac{\pa}{\pa z_{61}}, 
\vspace{0.2truecm}
\\
\zeta_3 & = & \frac{\pa}{\pa z_{32}} + z_{43}\frac{\pa}{\pa z_{42}} + \frac{1}{2}z_{43}^2\frac{\pa}{\pa z_{52}}, 
\end{array}
\right. 
\]
as a local frame of $E \subset T{\mathscr F}_{2,3}$. We set $E_1 = \langle \zeta_1\rangle, 
E_2= \langle \zeta_2, \zeta_3\rangle$. 

Then we have 
$[\zeta_1, \zeta_2] = 
\frac{\pa}{\pa z_{41}} + z_{43}\frac{\pa}{\pa z_{51}} - (z_{32}z_{43} - z_{42})\frac{\pa}{\pa z_{61}} =: \zeta_4, 
[\zeta_1, \zeta_3] = \frac{\pa}{\pa z_{42}} + z_{43}\frac{\pa}{\pa z_{52}} =: \zeta_5, 
[\zeta_2, \zeta_3] = 0$, that $[E_1, E_2] = E^{(2)} = \langle \zeta_i \mid 1 \leq i \leq 5\rangle$ and that 
$\rank(E^{(2)}) = 5$.  

Moreover we have 
$[\zeta_1, \zeta_4] = \frac{\pa}{\pa z_{51}} - z_{32}\frac{\pa}{\pa z_{61}} =: \zeta_6, 
[\zeta_1, \zeta_5] = \frac{\pa}{\pa z_{52}} =: \zeta_7, 
[\zeta_2, \zeta_4] = 0, [\zeta_2, \zeta_5] = 0, 
[\zeta_3, \zeta_4] = 0, [\zeta_3, \zeta_5] = 0$, so that 
$[E_1, E^{(2)}] = E^{(3)} = \langle \zeta_i \mid 1 \leq i \leq 7\rangle$ with $\rank(E^{(3)}) = 7$, 
while $[E_2, E^{(2)}] = E^{(2)}$. 

Further we have 
$[\zeta_1, \zeta_6] = 0, [\zeta_1, \zeta_7] = 0, 
[\zeta_2, \zeta_6] = 0, [\zeta_2, \zeta_7] = -\frac{\pa}{\pa z_{61}} =: \zeta_8, 
[\zeta_3, \zeta_6] = -\frac{\pa}{\pa z_{61}} = \zeta_8, [\zeta_3, \zeta_7] = 0$, 
that $[E_1, E^{(3)}] = E^{(3)}$, $[E_2, E^{(3)}] = E^{(4)}$, 
and that $\rank(E^{(4)}) = 8$, i.e. $E^{(4)} = T{\mathscr F}_{2,3}$. 

Therefore we see $E$ has $(3, 5, 7, 8)$ as the small growth vector. 
Moreover $E_1 = \langle \zeta_1\rangle, E_2 = \langle \zeta_2, \zeta_3\rangle$ 
gives a pseudo-product structure of type $B_3$ 
on $Z = {\mathscr F}_{2,3}$ (Definition \ref{B_3-pseudo-product}). 

\

(iii) (The $B_3(1,3)$-model) \ 
We study on $S = {\mathscr F}_{1,3}$. 
Let $V_1 = \langle h_1\rangle, V_3 = \langle h_1, h_2, h_3\rangle$, where 
\[
\left\{ 
\begin{array}{rcrrrrrrl}
h_1 & = & e_1 & + s_{21}e_2 & + s_{31}e_3 & + s_{41}e_4 + s_{51}e_5 + s_{61}e_6 + s_{71}e_7, 
\\
h_2 & = &      & e_2 &   & + s_{42}e_4 + s_{52}e_5 + s_{62}e_6 + s_{72}e_7, 
\\
h_3 & = &     &     &  e_3 & + s_{43}e_4 + s_{53}e_5 + s_{63}e_6 + s_{73}e_7. 
\end{array}
\right.
\]
The condition that $(V_1, V_3)$ is a null flag is given by
\[
\begin{array}{rcl}
(h_1 \vert h_1) & = & s_{71} + s_{21}s_{61} + s_{31}s_{51} - \frac{1}{2}s_{41}^2 = 0, 
\vspace{0.1truecm}
\\
(h_1 \vert h_2) & = & \frac{1}{2}s_{72} + \frac{1}{2}s_{21}s_{62} + \frac{1}{2}s_{31}s_{52} 
- \frac{1}{2}z_{41}z_{42} + \frac{1}{2}z_{61} = 0, 
\vspace{0.1truecm}
\\
(h_1 \vert h_3) & = & \frac{1}{2}s_{73} + \frac{1}{2}s_{21}s_{63} + \frac{1}{2}s_{31}s_{53} - \frac{1}{2}s_{41}s_{43} + \frac{1}{2}s_{51} = 0, 
\vspace{0.1truecm}
\\
(h_2 \vert h_2) & = & s_{62} - \frac{1}{2}s_{42}^2 = 0, 
\vspace{0.1truecm}
\\
(h_2 \vert h_3) & = & \frac{1}{2}s_{63} - \frac{1}{2}s_{42}s_{43} + \frac{1}{2}s_{52} = 0, 
\vspace{0.1truecm}
\\
(h_3 \vert h_3) & = & s_{53} - \frac{1}{2}s_{43}^2 = 0. 
\end{array}
\]
In particular $s_{21}, s_{31}, s_{41}, s_{51}, s_{61}, s_{42}, s_{52}$ and $s_{43}$ form a local system of coordinates of the null flag manifold ${\mathscr F}_{1,3}$. 

For a curve $(V_1(t), V_3(t))$, 
$V_1(t) = \langle h_1(t)\rangle, V_3(t) = \langle h_2(t), h_3(t)\rangle$ starting from $(V_1, V_3)$, 
$t \in (\R, 0)$, 
the distribution $S \subset T{\mathscr F}_{1,3}$ is defined by 
the condition that $V_1(t)' \subset V_2(t), V_2'(t) \subset V_2(t)^\perp$ and the condition 
is reduced that $h_1' \in V_2$ and only $(h_2' \vert h_3) = 0$. 
Then the condition is written in coordinates $s_{21}, s_{31}, s_{41}, s_{51}, s_{61}, s_{42}, s_{52}, s_{43}$ 
that 
\[
\begin{array}{c}
ds_{41} - s_{42}ds_{21} - s_{43}ds_{31} = 0, \quad ds_{51} - s_{52}ds_{21} - \frac{1}{2}s_{43}^2ds_{31} = 0, 
\vspace{0.1truecm}
\\
ds_{61} - \frac{1}{2}s_{42}^2ds_{21} + (s_{52} - s_{42}s_{43})ds_{31} = 0, \quad ds_{52} - s_{43}ds_{42} = 0. 
\end{array}
\]
Then we have a local frame of $L$  is given by 
\[
\left\{
\begin{array}{rcl}
\ell_1 & = & \frac{\pa}{\pa s_{42}} + s_{43}\frac{\pa}{\pa s_{52}}, 
\vspace{0.1truecm}
\\
\ell_2 & = & \frac{\pa}{\pa s_{43}}, 
\vspace{0.1truecm}
\\
\ell_3 & = & \frac{\pa}{\pa s_{21}} + s_{42}\frac{\pa}{\pa s_{41}} + s_{52}\frac{\pa}{\pa s_{51}} 
+ \frac{1}{2}s_{42}^2\frac{\pa}{\pa s_{61}}, 
\vspace{0.1truecm}
\\
\ell_4 & = & \frac{\pa}{\pa s_{31}} + s_{43}\frac{\pa}{\pa s_{41}} + \frac{1}{2}s_{43}^2\frac{\pa}{\pa s_{51}} 
- (s_{52} - s_{42}s_{43})\frac{\pa}{\pa s_{61}}. 
\end{array}
\right.
\]

Then $[\ell_1, \ell_2] = - \frac{\pa}{\pa s_{52}} =: \ell_5, 
[\ell_1, \ell_3] = \frac{\pa}{\pa s_{41}} + s_{43}\frac{\pa}{\pa s_{51}} + s_{42}\frac{\pa}{\pa s_{61}} =: \ell_6, 
[\ell_1, \ell_4] = s_{43}\frac{\pa}{\pa s_{61}} - s_{43}\frac{\pa}{\pa s_{61}} = 0, 
[\ell_2, \ell_3] = 0, 
[\ell_2, \ell_4] = \ell_6, 
[\ell_3, \ell_4] = 0$.  Therefore we have $L^{(2)} = \langle \ell_i \mid 1 \leq i \leq 6\rangle$ and 
$\rank(L^{(2)}) = 6$. 
Moreover we have 
$[\ell_1, \ell_5] = 0, 
[\ell_1, \ell_6] = \frac{\pa}{\pa s_{61}} =: \ell_7, 
[\ell_2, \ell_5] = 0, 
[\ell_2, \ell_6] = \frac{\pa}{\pa s_{51}} = : \ell_8
[\ell_3, \ell_5] = \ell_8, 
[\ell_3, \ell_6] = 0, 
[\ell_4, \ell_5] = - \ell_7, 
[\ell_4, \ell_5] = 0$. Thus we have $L^{(3)} = \langle \ell_i \mid 1 \leq i \leq 8\rangle = TS$ and 
$L$ is a $(4, 6, 8)$-distribution. 

We set $L = L_1\oplus L_2$, $L_1 = \langle \ell_1, \ell_2\rangle, L_2 = \langle \ell_3, \ell_3\rangle$, 
is a pseudo-product structure of type $B_3(1,3)$ in a generalised sense defined by Definition \ref{pseudo-product-B_3(1,3)}. 

\

(iv) (The $B_3(3)$-model) \ 
Lastly we consider $X = {\mathscr F}_3$ consisting of null $3$-spaces in $\R^{3,4}$. 
Set $V_3 = \langle k_1, k_2, k_3\rangle$ with 
\[
\left\{ 
\begin{array}{rcrrrrrrl}
k_1 & = & 
e_1 &  &  & + x_{41}e_4 + x_{51}e_5 + x_{61}e_6 + x_{71}e_7, 
\\
k_2 & = & 
& e_2 &  & + x_{42}e_4 + x_{52}e_5 + x_{62}e_6 + x_{72}e_7, 
\\
k_3 & = & 
 &     &  e_3 & + x_{43}e_4 + x_{53}e_5 + x_{63}e_6 + x_{73}e_7. 
\end{array}
\right.
\]
From that $V_3$ is a null $3$-space, we have 
\[
\begin{array}{rcl}
(k_1\mid k_1) & = & x_{71}  - \frac{1}{2}x_{41}^2 = 0, 
\vspace{0.1truecm}
\\
(k_1\mid k_2) & = & \frac{1}{2}x_{72} - \frac{1}{2}x_{41}x_{42} + \frac{1}{2}x_{61} = 0, 
\vspace{0.1truecm}
\\
(k_1\mid k_3) & = & \frac{1}{2}x_{73} - \frac{1}{2}x_{41}x_{43} + \frac{1}{2}x_{51} = 0, 
\vspace{0.1truecm}
\\
(k_2\mid k_2) & = & x_{62} - \frac{1}{2}x_{42}^2 = 0, 
\vspace{0.1truecm}
\\
(k_2\mid k_3) & = & \frac{1}{2}x_{63}  - \frac{1}{2}x_{42}x_{43} + \frac{1}{2}x_{52} = 0, 
\vspace{0.1truecm}
\\
(k_3\mid k_3) & = & x_{53} - \frac{1}{2}x_{43}^2 = 0. 
\end{array}
\]
In particular $x_{41}, x_{51}, x_{61}, x_{42}, x_{52}$ and $x_{43}$ form a local system of coordinates of the null Grassmannian ${\mathscr F}_3$. 
We define $D \subset T{\mathscr F}_3$ by the condition 
that $v \in D_{V_3} \subset T_{V_3}{\mathscr F}_3$ if and only if there is a curve $V_3(t)$ in ${\mathscr F}_3$ satisfying $V_3(0) = V_3$ and 
$V_3'(0) \in V_3^\perp$. The condition is written, by taking a frame $k_1(t), k_2(t), k_3(t)$ of $V_3(t)$, as 
\[
(k_1' \vert k_2) = 0, \ \ (k_1' \vert k_3) = 0, \ \  (k_2' \vert k_3) = 0, 
\]
at $t = 0$, that is, 
\[
-x_{41}'x_{42} + x_{61}' = 0, \ \ - x_{41}'x_{43} + x_{51}' = 0, \ \ - x_{42}'x_{43} + x_{52}' = 0, 
\]
at $t = 0$. 
Then we have that $D$ is given by the Pfaff system 
\[
dx_{61} - x_{42}dx_{41} = 0, \ dx_{51} - x_{43}dx_{41} = 0, \ dx_{52} - x_{43}dx_{42} = 0, 
\]
and obtain a local frame 
\[
{\textstyle 
\xi_1 = \frac{\pa}{\pa x_{41}} + x_{43}\frac{\pa}{\pa x_{51}} + x_{42}\frac{\pa}{\pa x_{61}}, \ \ 
\xi_2 = \frac{\pa}{\pa x_{42}} + x_{43}\frac{\pa}{\pa x_{52}}, \ \ \xi_3 = \frac{\pa}{\pa x_{43}}. }
\]
of $D$, satisfying bracket relations
\[
{\textstyle 
{[\xi_1, \xi_2]} = - \frac{\pa}{\pa x_{61}} =: \xi_4, \ \ 
[\xi_1, \xi_3] = - \frac{\pa}{\pa x_{51}} =: \xi_5, \ \ [\xi_2, \xi_3] =  - \frac{\pa}{\pa x_{52}} =: \xi_6, 
}
\]
with ${[\xi_i, \xi_j]} = 0 \ \ (0 \leq i \leq 3, 4 \leq j \leq 6)$. In particular $D \subset T{\mathscr F}_3$ is 
a $(3, 6)$-distribution.

\

\begin{flushleft}
Goo ISHIKAWA, 
\\
e-mail: ishikawa@math.sci.hokudai.ac.jp
\end{flushleft}

\begin{flushleft}
Yoshinori MACHIDA,
\\
e-mail: machida.yoshinori.a@shizuoka.ac.jp
\end{flushleft}

\end{document}